\newif\ifprivate                
\newif\ifconf                   
\newif\ifshort                  
\newif\ifcompress               
\newif\ifuc                     
\newif\iffinal                  

\privatefalse\conffalse\shortfalse\compressfalse\ucfalse\finalfalse

\documentclass[onecolumn]{revtex4-2}
\usepackage{tikz}
\usetikzlibrary{arrows.meta,positioning,decorations.pathreplacing,calc}

\newdimen\paravsp  \paravsp=1.3ex 
\usepackage[margin=1.5in]{geometry}
\usepackage{graphicx}
\usepackage{epstopdf}
\usepackage{amsmath}
\usepackage{amsfonts}
\usepackage{bbm}
\usepackage{braket, tabularx}
\usepackage{rotating}  

\usepackage{color}

\usepackage{lineno}

\usepackage{comment}

\usepackage[normalem]{ulem}

\usepackage{color}
\usepackage{subfigure}

\makeatletter

\newcommand{\Rmnum}[1]{\expandafter\@slowromancap\romannumeral #1@}
\makeatother

\ifprivate
\newcommand{\kamal}[1]{{\color{magenta}  #1}}
\newcommand{\ard}[1]{{\color{red}  #1}}
\else
\newcommand{\kamal}[1]{}
\newcommand{\ard}[1]{}
\fi

\def\eps{\varepsilon}            
\def\paradot#1{\vspace{\paravsp plus 0.5\paravsp minus 0.5\paravsp}\noindent{\bf\boldmath{#1.}}}

\def\fr#1#2{{\textstyle\frac{#1}{#2}}} 
\def\frs#1#2{{^{#1}\!/\!_{#2}}} 

\ifprivate\def\todo#1{{\color[rgb]{.6,.3,.0}\it[[ToDo:\ #1]]}}\else\def\todo#1{}\fi 
\ifprivate\def\note#1{{\color[rgb]{.5,.5,.0}\it[[Note:\ #1]]}}\else\def\note#1{}\fi 

\def\cX{\mathcal{X}}
\def\cN{\mathcal{N}}

\begin{document}
\title{Retrodicting Chaotic Systems: An Algorithmic Information Theory Approach}

\author{Kamal Dingle$^{1}$, Boumediene Hamzi$^{2,3,4}$, Marcus Hutter$^{5,6}$,  Houman Owhadi$^{2}$}

\affiliation{
$^1$Centre for Applied Mathematics and Bioinformatics,
Department of Mathematics and Natural Sciences, 
Gulf University for Science and Technology, Kuwait\\
$^{2}$Computing and Mathematical Sciences Department, Caltech, USA \\
$^{3}$The Alan Turing Institute, London, UK \\
$^4$Isaac Newton Institute for Mathematical Sciences, Cambridge, UK.\\
$^{5}$Google DeepMind, London, UK\\
$^{6}$School of Computing, Australian National University, Canberra, ACT 2601, Australia 
}

\date{\today}

\begin{abstract}
\noindent 
Making accurate inferences about data is a key task in science and mathematics. Here we study the problem of \emph{retrodiction}, inferring past values of a series, in the context of chaotic dynamical systems. Specifically, we are interested in inferring the starting value $x_0$ in the series $x_0,x_1,x_2,\dots,x_n$ given the value of $x_n$, and the associated function $f$ which determines the series as $f(x_i)=x_{i+1}$.  Even in the deterministic case this is a challenging problem, due to mixing and the typically exponentially many  candidate past values in the pre-image of any given value $x_n$ (e.g., a current observation). We study this task from the perspective of algorithmic information theory, which motivates two approaches: One to search for the `simplest' value in the set of candidates, and one to look for the value in the lowest density region of the candidates. We test these methods numerically on the logistic map, Tent map, Bernoulli map, and Julia/Mandelbrot map, which are well-studied maps in chaos theory. The methods aid in retrodiction by assigning low ranks to candidates which are more likely to be the true starting value. Our approach works well in some parameter and map cases, and outperforms several other retrodiction techniques (each of which fails to outperform random guessing). Nonetheless, the approach is not effective in all cases, and several open problems remain including computational cost and sensitivity to noise. All of these methods are unified through a Gaussian Process (GP) perspective, motivating complexity-based priors for GPs. \\


\end{abstract}

\maketitle




\section{Introduction}\label{sec:intro}

A core task in science, mathematics, and indeed more broadly, is generalizing patterns in data to infer unknown or unseen values. Examples include predicting the outcome of a physics experiment, or inferring the mechanisms that caused some animal species to die out. 
Much of time series analysis is concerned with developing strategies and methods for forecasting, which is to say, predicting the future values of a series. The related problem of \emph{retrodiction} --- inferring past values from present observations --- has received perhaps less attention, but nonetheless has been studied for decades \cite{gibbs1902elementary,watanabe1955symmetry,ellison2009prediction,rupprecht2018limits,surace2022state,parzygnat2023axioms,crutchfield2009time}. 

Prediction and retrodiction, although both viewed as inferring a missing value from  a (time) series, are not necessarily equivalent problems with similar methods and challenges,  or necessarily equally easy. In some senses, from a theoretical viewpoint, retrodiction can be more difficult than prediction in fact, as Rupprecht and Vural show \cite{rupprecht2019enhancing}.  On the other hand, 
Crutchfield et al. \cite{crutchfield2009time} show that prediction and retrodiction are equally easy for a stationary process, but that predicting and retrodicting may require different amounts of information storage.
Interestingly, for humans, there are reports that find prediction and retrodiction are of equal difficulty \cite{jones2007mind,smith2014looking}, and also that participants are better at inferring unseen past (versus future) events \cite{xu2024temporal}, and that humans employ a more complex mix of strategies in certain environments \cite{sharp2024humans}.

In this work, we study the problem of retrodiction in time series obtained from deterministic chaotic dynamical systems.  More specifically, we consider the series $x_0,...,x_n$ generated by $n$ iterations of some given map $f$, with starting value $x_0$. The task is to infer $x_0$ from $x_n$. Our approach could also be called \emph{exact} retrodiction, because we try to identify the exact starting value $x_0$ from a finite list of possibilities, rather than estimating the approximate value of $x_0$ via fitting to trends, or using coarse-grain methods like binning which make the dynamics stochastic.
This retrodiction problem is challenging because the concept of mixing from ergodic theory implies that any given currently observed value can come from almost anywhere \cite{hasselblatt2003first}; in other words the pre-images of a given point $x_n$ will typically be well distributed on the phase space, and therefore we have little knowledge about where the given point originated from. Moreover, as we will show, the finite list of possible starting values is typically exponentially large, exacerbating the problem of identifying the true starting value $x_0$.

While similar problems have been studied earlier \cite{rupprecht2018limits} (and citations above), our approach is quite different from such works, because we frame the problem in terms of \emph{algorithmic information theory} (AIT) \cite{solomonoff1960preliminary,solomonoff1964formal,kolmogorov1965three,chaitin1975theory}. AIT is a subfield of theoretical computer science, and concerns the information content of individual discrete data objects, like binary strings, integers, graphs, etc. By viewing the problem of retrodiction through the lens of AIT and \emph{Kolmogorov complexity}, we develop approaches to the (exact) retrodiction problem described above. We apply our methods to numerical simulations of the logistic map, Tent map, Bernoulli map, and Julia/Mandelbrot map. In some map/parameter cases the methods are effective in the sense of assigning low rank values to possible starting values, thereby aiding in retrodiction. In some other cases however, the methods do not help to reduce the size of the possible candidates, and reasons for this are discussed. We also compare to several retrodiction methods such as Gaussian Processes, and others, each of which fails to perform better than random guessing. This failure highlights the difficulty of the current problem, and the superiority of our proposed methods under some circumstances.

Another contribution of this work  is the unification of the various retrodiction methods through a common Gaussian Process (GP) framework. By leveraging complexity-based priors for GPs, we introduce a unifying perspective that bridges machine learning with algorithmic information theory. The connection between AIT and mainstream machine learning, particularly through GP models with complexity-driven priors, offers a novel direction for further research, as explored in recent papers by Hamzi, Hutter, and Owhadi \cite{Hamzi2024, Hamzi2024b, Hamzi2025, HamziHutterOwhadi2025_SGP}


 Overall, this work contributes to solving this challenging problem of retrodiction in the context of chaotic dynamics \cite{pellicer2015importance}, but does not claim to have settled the question.

\section{Problem set-up}\label{sec:prob}

We will now introduce notation and describe the proposed problem in more detail.

\paradot{Definitions}
Consider a surjective map $f:\cX\to\cX$, which generates sequence 
\begin{equation}
    x\equiv x_0,...,x_n\equiv y
\end{equation}
from $x$ via $x_{i+1}=f(x_i)$.
Typically $\cX\subseteq\mathbb{R}^d$ for some $d\in\mathbb{N}$. The central challenge is this: Given $y$ we want to infer $x$. If $f$ is injective and known, then this is in principle possible:
Let 
\begin{equation}
    y ~=~ g(x) ~=~ f(f(...f(x)))
\end{equation}
be the $n$-times iteration of $f$, then
\begin{equation}
    x ~=~ \bar g(y) ~:=~ g^{-1}(y)
\end{equation}
If $f$ is not injective, then $f^{-1}$ and hence $g^{-1}$ are not unique. 
In this case, let 
\begin{equation}
    \cX_0 ~:=~ \bar g(y) ~=~ \{x^1,...,x^m\}
\end{equation}
be all pre-images of $g$, and $m\in\mathbb{N}$ is the number of elements in $\cX_0$. (Note that $x^i$ should not be confused with the $i$th iterate $x_i$). We will refer to $\cX_0$ as the \emph{candidate} values for $x_0$.
 If $f$ always has two pre-images, then $|\cX_0|=m=2^n$. 

Because the true starting value $x\in\cX_0$, so given $y$ and $f$,
we can determine $x$ up to $\log_2 m$ bits. In other words, we can quantify our uncertainty in the value of $x$ by equating complete ignorance of which element of $\cX_0$ is $x$ with having a uniform distribution over the $m$ elements; and this uniform distribution would have $\log_2 m$ bits of Shannon entropy/uncertainty. The question we wish to study in this work is, can we do better than $\log_2 m$ bits of uncertainty? Naturally, if we can certainly identify the correct value of $x$ then this would reduce our uncertainty to zero bits. So we aim to develop methods that can lower our uncertainty from $\log_2 m$ bits to as close to zero as possible.

The problem we study could also be called \emph{exact} retrodiction, because we try to identify the exact starting value $x_0$ from a (finite) collection of $m$ elements in $\cX_0$, rather than estimating the approximate value of $x_0$, as is often done in time series analysis \cite{chatfield2019analysis}.

There are many possible ways to approach the problem of retrodiction, even if we restrict to deterministic dynamical systems. The approach we will propose here is based on assigning a rank to each element of $\cX_0$ according to how likely it is to be the starting value $x_0$. Then we can encode $x_0$ by its rank $r\in\mathbb{N}$ using $\log_2 r$ bits. If $\log_2 r \ll \log_2 m$, then many bits would have been saved and we have much less uncertainty about the starting value than $\log_2 m$ bits. On the other hand, if $\log_2 r \approx \log_2 m$, then barely any bits would be saved, and our ignorance about the starting value would not have decreased.

\paradot{Comment:} 
Note that $\log r$ and $\log_2 m$ are not the lengths of prefix-free codes \cite{cover2006elements}, so we need to be careful in comparing similar types of code lengths. We can either work without prefix codes, and just compare $\log_2 r$ to $\log_2 m$, or we could write the rank as a prefix code as follows: We can get a prefix code length by first assigning a probability $Q(r)$ to $r$, and then take $CL(r):=\log_2 Q(r)^{-1}$ as code length.
The easiest way is to construct the CDF $F$ or more conveniently the survival function $\bar F=1-F$ of $Q$,
then $Q(r)=\bar F(r-1)-\bar F(r)$ for $r\in\{1,...,m\}$.
We would require $\bar F$ to be monotone decreasing from $\bar F(0)\leq 1$ and $\bar F(m)\geq 0$,
with $=$ giving a complete code.
$F(r)=\fr{r}{m}$ gives $Q(r)=\fr1m$, hence the naive uniform CL$(r)=\log_2 m$.
For $\bar F(r)=\fr1{r+1}$ we get $Q(r)=\fr1{r(r+1)}$, 
hence CL$(r)=\log_2 r(r+1)$ which is about double the plain code length $\log_2 r$ 
and not good enough for our purpose.
$$
  \bar F(r) ~:=~ h(r)/h(0) ~~~\text{with}~~~ h(r):=[\ln(r+a)]^{-1}-[\ln(m+a)]^{-1}
$$
for some $a>1$ gets rid of the factor $2$ for the price of an additive loglog term:
$$
  \text{CL}(r) ~\approx~ \log_2(r+a) + 2\log_2\ln(r+a) - \log_2\ln(a) 
  ~\approx~  \log_2 r + 2\log_2\ln r
$$
If we have some prior knowledge, 
then $a\approx<r>/\ln<r>$ is a good choice, where $<r>$ is the average $r$.

Having made this discussion of prefix-codes, for the purpose of this present work, we will work with plain codes and just compare $\log_2 r$ to $\log_2 m$, where the rank is $1\leq r\leq m$.

\paradot{Comment}
While we will not use this approach in this current work, we point out for potential future extensions that a different way to assign a code length is to assign a probability $p_x=P(x|y)$ to $x\in\cX_0$.
The probabilities may arise from some unnormalized `density' $\rho(x)$ 
as $p_x:=\rho(x)/\sum_{x\in\cX_0}\rho(x)$.
We can then compare the naive probability $\fr1m$ that $y$ arose from $x_0$ to $p_{x_0}$.
Alternatively, on a log-scale we can check whether the Huffman or Shannon-Fano code of $x_0$ w.r.t.\ $p_x$
of length $\log_2 p_x^{-1}$ is smaller than a uniform code of length $\log_2|\cX_0|=\log_2 m$.
Strictly speaking the code has length $\lceil\log_2 p_x^{-1}\rceil\in\mathbb{N}$,
but if we jointly encode multiple items, we only lose at most $1$ bit overall,
which justifies dropping $\lceil\rceil$.

\section{Maps used in the study}

\subsection{Maps and example}
We will study (exact) retrodiction in the following well-known maps from dynamical systems:
The Logistic map \cite{lorenz1964problem} is perhaps the most iconic map in dynamical systems, and has been extensively studied \cite{berger2001chaos}. It is a 1D map, and has the form 
\begin{equation}
    f(x) = \mu x (1-x)
\end{equation}
with $x\in[0,1]$ and $\mu \in [0,4]$, and therefore $f^{-1}(x)=0.5 \pm \sqrt{0.25 - x/\mu}$.
We will use $\mu=4.0$ in this study. Notice that for this map there are typically two distinct possible values of $x_{n-1}$, and four distinct possible values of $x_{n-2}$, continuing to $2^n$ possible values of $x_0$. It is well known that the limiting distribution of values for $x_n$ is the Beta(0.5,0.5) distribution, which has the density $1/(\pi\sqrt{x(1-x)})$. This means that for large $n$, most of the iteration values $x_n$ are close to 0 or 1, with a low probability of being close to 0.5. This distribution also applies for backward iterations, so that for typical values in [0,1], the set of candidates in the pre-image will also follow the Beta(0.5,0.5) distribution.

The Tent map \cite{hasselblatt2003first} is another commonly studied 1D dynamical system. It has the form
\begin{equation}
     f(x) = \begin{cases} 
          2x, & 0.0\leq x\leq 0.5\\
          2-2x, & 0.5\leq x\leq 1.0\\
       \end{cases}
\end{equation}
with $x\in[0,1]$, and the pre-images can be found with $f^{-1}(x)= x/2, 1-x/2$. It is well known that the limiting distribution of values for $x_n$ is the uniform distribution over [0,1]. This also applies for backward iterations, so that for typical values in [0,1], the set of candidates in the pre-image will also follow the uniform distribution over [0,1]. This means that in terms of retrodiction, we have essentially no idea which intervals are more or less likely to contain the true starting value $x_0$ (assuming the number of iterations is large enough).

The Bernoulli map \cite{hasselblatt2003first} shares many similarities with the Tent map, and has the form 
\begin{equation}
f(x)= (2x \mod 1)    
\end{equation}
with $x\in[0,1]$, and the pre-images are found by $f^{-1}(x)= x/2, 0.5 + x/2$. It is well known that the limiting distribution of values for $x_n$ is the uniform distribution over [0,1]. Like the Tent map, for typical values in [0,1], the set of candidates in the pre-image will also follow the uniform distribution over [0,1]. This property makes retrodiction very challenging.

Finally, the Julia map, or Mandelbrot quadratic map, is defined by
\begin{equation}
f(z)=z^2 +c    
\end{equation}
so that the two pre-images are given by $f^{-1}(z)=\pm \sqrt{z-c}$ and $z,c\in \mathbb{C}$. Because $\mathbb{C}$ is geometrically equivalent to $\mathbb{R}^2$, the Julia map can be described as a 2D map (as opposed to the other maps which are 1D). The limiting distribution of the values of iterates of this map depends on the values of $z_0$ and $c$, and hence, there is no one general distribution that is applicable to all relevant values of $z,c\in \mathbb{C}$.

Note that in all these maps there are typically two distinct values in the pre-image of $x$ (or $z$). Therefore we typically have $m=2^n$ candidate starting values. The exponentially large set of candidate values makes the problem of exact retrodiction hard.

{\bf Example:} To fix thoughts, we give a concrete example of the retrodiction problem we consider in this article. Suppose $x=x_0=0.35$ is the (unknown) starting value in the logistic map with $\mu=4$, and the number of iterations is $n=3$. This gives the series 
\[x_0=0.35, x_1=0.91, x_2=0.3276,  x_3=0.88111296\]
The problem is this: We observe $y=x_3$ as the  current value of the series and we want to perform retrodiction of the series to infer the value of $x_0$. This is a nontrivial problem because taking the pre-images of $x_3=0.88111296$, we get the following list of $m=2^n=2^3=8$ candidate values: 
\[ \cX_0 = \{0.7312..., 0.2688..., 0.94333..., 0.05667..., 0.65, 0.35, 0.97697..., 0.02303...\} \]
Because we know the map, and it is assumed to be deterministic, the true (unknown) starting value $x_0=0.35$ is indeed in this set. The problem we study is how to identify the true $x_0$ value from the (exponentially large) list of candidates in $\cX_0$. 
 Even if we cannot identify $x_0$ itself, we attempt at least to reduce the uncertainty in from $\log_2(m)$ bits, which would be the case if a uniformly random sampling search method was employed. Hence $\log_2(m)$ bits acts as a baseline to improve on.

\subsection{How $x_0$ is chosen}
The manner in which $x_0$ is initially chosen is important, and affects our ability to infer (retrodict) its value. We will consider several possibilities. 

If we knew that $x_0$ was chosen by a human, then we might have some prior knowledge about the kinds of values that are typically selected (or not) by humans. In principle, this could help with the retrodiction problem, by for example invoking some psychological knowledge. However, we will not explore this further as it is not of primary interest.
If $x_0$ is chosen by first iterating the relevant function $f$ many times on some arbitrary choice of starting value, then we would have prior knowledge about the likely value of $x_0$, based on knowing the typical distribution of iterations of that map. For example, if it was the logistic map, then we could guess that the value might follow the Beta$(0.5,0.5)$ distribution.
Another conceivable way for $x_0$ to be chosen is to start with $x_n$, and then one of the candidates $\cX_0$ is chosen at random to define $x_0$. This would be an impossible problem to solve, because there would be no way to know which is the true value of $x_0$, beyond simply guessing at random uniformly from the candidates.  We are not concerned with this manner of choosing $x_0$, but instead will assume $x_0$ is a fixed finite complexity value, while $n$ can vary and grow large.

\section{Approach 1: complexity of candidates}\label{sec:meth1}

Having defined the problem, and introduced the maps we will study, we will describe the retrodiction methods we introduce. In order to do so, we must first give some background details on the relevant theory. 

\subsection{Background on Kolmogorov complexity}

Developed within theoretical computer science by Solomonoff, Chaitin, and Kolmogorov, \emph{algorithmic information theory} \cite{solomonoff1960preliminary,kolmogorov1965three,chaitin1975theory} (AIT)  connects computation, computability theory and information theory. The central quantity of AIT is \emph{Kolmogorov complexity}, $K(x)$, which measures the complexity of an individual object $x$ as the amount of information required to describe or generate $x$.  More formally, the Kolmogorov complexity $K_U(x)$ of a string $x$ with respect to a universal Turing machine \cite{turing1936computable} (UTM) $U$,  is defined \cite{solomonoff1960preliminary,kolmogorov1965three,chaitin1975theory} as
\begin{equation}
K_U(x) = \min_{p}\{|p|: U(p)=x\}
\end{equation}
where $p$ is a binary program for a prefix (optimal) UTM $U$, and $|p|$ indicates the length of the (halting) program $p$ in bits.  Due to the invariance theorem \cite{li2008introduction} for any two optimal UTMs $U$ and $V$, $K_U(x) = K_V(x)+O(1)$ so that the complexity of $x$ is independent of the machine, up to additive constants. Hence we conventionally drop the subscript $U$ in $K_U(x)$, and speak of `the' Kolmogorov complexity $K(x)$. Despite being a fairly intuitive quantity and fundamentally just a data compression measure, $K(x)$ is uncomputable, meaning that there  cannot exist a general algorithm that for any arbitrary string returns the value of $K(x)$.
Informally, $K(x)$ can be defined as the length of a shortest program that produces $x$, or simply as the size in bits of the compressed version of $x$. If $x$ contains repeating patterns like $x=1010101010101010$ then it is easy to compress, and hence $K(x)$ will be small. On the other hand, a randomly generated bit string of length $n$ is highly unlikely to contain any significant patterns, and hence can only be described via specifying each bit separately without any compression, so that $K(x)\approx n$ bits. $K(x)$ is also known as  \emph{descriptional complexity}, \emph{algorithmic complexity}, and \emph{program-size complexity}, each of which highlight the idea that $K(x)$ measures the amount of information required to describe or generate $x$ precisely and unambiguously. 
More details and technicalities can be found in standard AIT references \cite{li2008introduction,calude2002information,gacs1988lecture,shen2022kolmogorov}.

Directly applying AIT and algorithmic probability to practical problems such as time series is not straightforward due to a number of reasons. These include the fact that Kolmogorov complexity is uncomputable, the theory is asymptotic with equations given only up to $O(1)$ terms, and the theory is developed in terms of universal Turing machines \cite{turing1936computable}, which are not common in practical settings. Despite these and other issues, there are arguments that motivate the application of AIT in practical settings in combination with various approximations and assumptions, and very many studies have applied AIT, including to make useful and verified predictions. Some examples are in physics  \cite{bennett1982thermodynamics,kolchinsky2020thermodynamic,zurek1989algorithmic,mueller2020law,avinery2019universal,martiniani2019quantifying}, biology   \cite{ferragina2007compression,adams2017formal,johnston2022symmetry,dingle2022predicting}, networks \cite{zenil2014correlation,zenil2018review}, in addition to data analysis \cite{vitanyi2013similarity,cilibrasi2005clustering,zenil2019causal,zenil2011algorithmic,dingle2023note}, among several others \cite{li2008introduction,dingle2018input,dingle2020generic}.

\subsection{Kolmogorov complexity of candidates}

Our first approach to the retrodiction problem comes from examining the Kolmogorov complexity of candidates.
In this section, we will assume that  $x$ is describable with finite complexity, i.e., $K(x)<\infty$. This would be the case if $x$ is specified to $D$ decimal places, where $D$ is finite, or for numbers like $\pi$ or $\sqrt{2}$ which have infinite decimal expansions, yet are finitely describable via simple formulas.  

At a high level, the idea of this method is to show that because the complexity of typical candidates grows with $n$, then for any starting value of finite complexity, eventually, the true starting value will be contained within a relatively small subset of low-complexity candidates. Thus while the set of candidates grows exponentially, the true starting value may be much easier to find than `looking for a needle in a haystack'. Now we give a more detailed explanation.

We can estimate the typical complexity of elements of $\cX_0$ as follows:
Any member $x^i\in\cX_0$ can be described by describing $y$, then using $f$ and $n$ to generate the set $\cX_0$, and then specifying the index $i$ of element  $x^i$ using $K(i)$ bits. It follows that we can upper-bound the complexity of $x^i$ as
\begin{equation}
    K(x^i) \leq K(f) + K(n) + K(y)+ K(i)+O(1)
\end{equation}
with $1\leq i \leq m$. The complexity of the function $f$ will be assumed to be low, because many dynamical systems (e.g., the Tent map) have simple maps, which take only a few bits to describe. Hence we can assume $K(f)=O(1)$. For typical $n$, we have $K(n)\approx \log_2(n)$ (ignoring loglog terms), which is typically much smaller than $K(i)$, and so could be ignored. It follows that  typical candidates have complexity 
\begin{equation}
    K(x^i) \lesssim   K(y) + K(m)\label{eq:y_m}
\end{equation}
because $K(i)\approx K(m)$ for almost all $i$. The complexity of $y$ grows slowly with $n$, and can be bounded as 
\begin{equation}
    K(y)\leq K(x_0) + K(n)+K(f)+O(1)
\end{equation}
Hence we can say that for relatively simple $y$, or alternatively for large enough iterations $n$, we can reduce the approximation in Eq.\ (\ref{eq:y_m}) to just 
\begin{equation} 
   K(x^i) \sim K(m) \sim  n \label{eq:m_n}
\end{equation}
if we assume that $|\cX_0|=m=2^n$. The bound in Eq.\ (\ref{eq:y_m}) has been replaced by an approximate equality in Eq.\ (\ref{eq:m_n}), because it is well known \cite{li2008introduction} that in a set of $k$ different elements, the complexity of typical members of the set must be at least $O(\log_2(k))$.
Thus we have shown that the complexity of typical members of $\cX_0$ grows roughly linearly with $n$ (given a few assumptions). By contrast, the specific starting value $x_0$ is (by assumption) some fixed complexity number. It follows that for $n$ large enough we have
\begin{equation}
    K(x_0)\ll K(x^i)
\end{equation}
for typical $x^i$; i.e., $x_0$ will be much simpler than the majority of other members of $\cX_0$. This also means that the required size of $n$ for the method to be effective for starting value $x_0$ is $K(x_0)\lesssim n$. 

It is possible that there are other  elements of $\cX_0$ which are simpler than $x_0$ itself, but such elements will be relatively rare. To see why, note that 
\begin{equation}
    |\{x^i: K(x^i)<K(x_0)\}|\lesssim 2^{K(x_0)}\label{eq:size_of_set}
\end{equation} 
because the number of objects of complexity $k$ must be less than $2^k$ \cite{li2008introduction}, so that the fraction of candidates which are simpler than $x_0$ must vanish as 
\begin{equation}
    \frac{2^{K(x_0)}}{m}\rightarrow 0
\end{equation}
as $n\rightarrow \infty$. Hence, the true starting value will be in a `small' set of low-complexity candidate values.

\subsection{Retrodiction approach}

The above arguments suggest an approach to the retrodiction problem: enumerate the elements of $\cX_0$, and order them from simplest to most complex. Those with lower rank are more likely to be the true starting value $x_0$. (If desired, one can sample the sorted candidates with probability taken inversely to the rank in the sorted list, thereby increasing the chances of correctly guessing the true starting value.)
Unfortunately, we cannot compute Kolmogorov complexity for each candidate, because it is formally uncomputable \cite{li2008introduction}. Nonetheless, as is by-now well established, in practical situations approximations to Kolmogorov complexity often work well \cite{li2008introduction,vitanyi2013similarity,dingle2018input,vitanyi2020incomputable}, which motivates finding some reasonable (computable) complexity estimator \cite{janzing2010causal} and testing the method numerically.

One way to estimate the complexity of a decimal number candidate in [0,1] is by using a standard lossless data compression technique, and compress the values of the candidates. Here we will use the \texttt{zlib.compress} package in Python, which employs essentially the same compression algorithm as in gzip (i.e., the `deflate' method which combines Lempel-Ziv 1977 and Huffman coding). Following this, if the compressed size for $x_0$ is less than the compressed size of the decimal expansion of typical members of $\cX_0$, then $x_0$ will stand out as being simple. 
Note that the method described above is not fundamentally based on one compression algorithm or another, but rather on the Kolmogorov complexity of candidates. Hence, building on this current work, in the future, other (more refined) approximations of complexity might be developed and applied using the same method described here.

Interestingly, while the above argument does not \emph{invoke} Occam's razor \cite{sober2015ockham}, it does point to an analogous conclusion: choosing the simplest candidate is a good strategy for inferring the true starting value. 

\subsection{Expected rank}

It is interesting to consider the expected value of the rank $r$ as a function of the complexity of $x_0$. This question has a couple of angles. One is theoretical, namely rank in terms of Kolmogorov complexity, and another concerns the accuracy of the complexity estimate. 

In the case that $x_0$ is a random real number, then it has infinite Kolmogorov complexity, $K(x_0)=\infty$. In this case, assuming the map yields $\approx 2^n$ distinct candidates and they all have infinite complexity, then comparing complexity values of the candidates will be meaningless and no retrodiction is possible. Hence the expected rank is
\begin{equation}
    \mathbb{E}[r]=(2^n+1)/2
\end{equation}
the same value as for uniform random guessing out of $2^n$ possibilities. 

As for the case when $n\approx 1$ is small and the complexity of $x_0$ is finite, $K(x_0)<\infty$, all the candidates will have similar complexity values due to the following argument:  If $n\approx 1$ then $K(x_0)\approx K(x^i)$ for any $i$ because any $x_i$ can be described using $x_0$ and iterating forward  to $x_n$, and then specifying any $x^i$ out of the (small) set of $2^n$ candidates. Because in this case all candidates have similar (finite) complexity, then $\mathbb{E}[r]\approx (2^{n}+1)/2$ and retrodiction success is no better than random guessing.

In the case that the complexity of $x_0$ is finite, $K(x_0)<\infty$, $n$ is large, and the map yields $\approx 2^n$ distinct candidates, then we have Eq.\ (\ref{eq:size_of_set}) to guide estimates: Within the $2^n$ candidates, at most $2^{K(x_0)}$ will have rank lower than $x_0$, hence 

\begin{equation}
\mathbb{E}[r]\lesssim 2^{K(x_0)}\ll 2^n    
\end{equation}
Whether or not the expectation for a given map will be much less than $2^{K(x_0)}$ depends on the distribution of complexities within $\cX_0$: It may be that there are many elements of $\cX_0$ that have a complexity value less than $K(x_0)$, in which case $\mathbb{E}[r]\approx 2^{K(x_0)}$ out of $2^n$ candidates. Alternatively, there may be no or very few elements in $\cX_0$ that are simpler than $K(x_0)$, in which case $\mathbb{E}[r]\approx 1$. A general prediction about which of these (or other) cases will appear depends on the details of a given map. 

The preceding analysis has assumed the use of the true Kolmogorov complexity, but in practice any implementation of this proposed method will make use of some real-world complexity estimator such as Lempel-Ziv compression. Such complexity measures introduce another level of uncertainty, namely their accuracy in assigning a correct complexity value. It may be, and indeed it is highly likely, that a complexity measure occasionally assigns a high complexity to an object which is, in fact, low complexity. Hence, even if in theory the starting value $x_0$ is low-complexity, and theoretically we should have the rank $r\ll 2^n$, then it may be that the rank is $r\approx 2^n$ due to limitations of the complexity measure. It is not straightforward to determine the frequency with which this will happen, not least because it depends on the specific complexity measure used.

\subsection{Limitations of approach}\label{limitations}

The proposed method has some limitations, which we discuss now.
To begin, the method requires enumerating an exponentially large set of candidates, and then searching all of these elements for simpler values. Hence it is computationally expensive for larger $n$. Having said this, we hold that performing retrodiction is interesting and important even for non-asymptotically large $n$, even for say $n\approx 15$. Hence, while this is a limitation of the method, it does not preclude its utility completely. Further, for future work it might be possible to propose ways to trim
 down the size of the set so as to reduce the computational cost, but we leave this for later investigation.

Another limitation is that the method is not very practical for the real world, e.g., in physics, where 100\% accuracy in decimal numbers is not plausible. Similarly, random perturbations can disrupt the method, or a series may be deterministic but noise may be introduced during measurement. Again, while this is a limitation, we hold that the approach is a worthwhile first step in solving the very challenging problem of retrodiction in chaotic systems. At present, we are not proposing to implement the method in e.g. physics and engineering problems, but study it as a novel mathematical investigation into viewing the problem from the perspective of AIT. If successful at this stage, it will motivate refinements and adjustments for future applications. It is plausible that with very small random noise, and with small $n$, the method could still be effective (see more on this below). 

Numerical issues can be a problem for this method also, in that when iterating maps forward and backward, numerical errors can be introduced, which essentially introduce a small amount of randomness into the candidate values. If these errors accumulate and become large, then they can cause the method to fail. One way to handle this is to use higher-precision computational tools, such as \texttt{Decimal} in Python, and set the precision level high.

Because Kolmogorov complexity is uncomputable, naturally this poses a challenge to a strategy using Kolmogorov complexity. However, as mentioned above, there are by now very many studies which use approximations to Kolmogorov complexity in practical settings with success \cite{li2008introduction,Li:07appait}. So, while better complexity measures would aid, we can still expect some success while using computable approximations.

It is well known that almost all real numbers have infinite Kolmogorov complexity \cite{li2008introduction}, because of their infinite random digits. Hence, for a typical random real value $x$, the complexity of $x$ and the complexity of $f(x)$ will be the same (i.e., infinite), and indeed presumably all the candidate values will also have infinite complexity, and hence there is no possibility of identifying a simpler candidate. On the one hand, this observation means that the method we propose would not work, in principle, for nearly all real numbers. On the other hand, this is not a significant problem because in practical computation, we typically work with truncated or finite-complexity numbers. Additionally, the problem of exact retrodiction even for finite-complexity values is still interesting and a worthwhile studying.

\subsection{An algorithmic probability argument}\label{sec:aparg}

We also propose another perspective on the retrodiction problem, also fundamentally based on complexity. We can embed the problem at hand into a more general question: Given only \emph{some} elements of a sequence $x_0,...,x_n$, how shall we reconstruct the missing values?
Let $U\subseteq\{0,...,n\}$ be the index set of unknown $x_i$.
If $U=\{n\}$, this is a classic time-series forecasting problem, which can be solved in the most general case via Solomonoff induction based on Kolmogorov complexity quantifying Occam's razor \cite{solomonoff1960preliminary,li2008introduction,hutter2004universal,hutter2007universal}.
We can try to apply Occam's razor also for general $U$ telling us to choose the simplest sequence consistent with what we know.
Quantified in terms of Kolmogorov complexity that would be the sequence  that achieves the minimal length 
\begin{equation}
    \ell^*:=\min_{\{x_i:i\in U\}} K(x_0...x_n)
\end{equation}

In our particular setup, $U=\{0,...,n-1\}$ and $K(x_0...x_n|f)=K(x_0|f)$,
or in reduced form, $K(x,y|g)=K(x|g)$, since $x_1...x_n$ can be computed from $x_0$.
Conditioning on $y$ we get $K(x,y|y,g)=K(x|y,g)+O(1)$, hence 
\begin{equation}
\ell^*=\min_x K(x|y,g)    
\end{equation}
The algorithmic probability formulation would lead to
\begin{equation}
    P(x|y)=2^{-K(x|y,g)}
\end{equation}
 This `derivation' has weaknesses and one has to be careful with $O(1)$ and sometimes even $O(\log)$-terms,
but it can serve as a motivation for choosing the simplest value, and the more practical approaches below. 
To expound on this, because AIT results typically involve $O(1)$ terms, and these terms are not possible to remove \cite{muller2010stationary}, making precise quantitative estimates can be challenging (unless asymptotically large variables are assumed). Despite this, there is a whole area of applied AIT which operates by using Kolmogorov complexity bounds or equations (which involve $O(1)$ terms), and then using these as guides for more practical implementations of the theory \cite{vitanyi2013similarity,cilibrasi2005clustering}. This approach has proved very successful, and we follow suit.
See also See \cite[Sec.8.2.2]{hutter2004universal} and \cite[Sec.5.9]{Hutter:11uiphil} for a more detailed discussion of the $O(1)$ terms.

The generic idea is to have a bias towards $x$ that have a simple description given $y$ and $g$.
In theory, Kolmogorov complexity is the ultimate description length,
in practice, we can choose descriptions more tailored to the problem at hand (as discussed above). 
Due to the conditional complexity term, we might also consider $x$ values that are similar to $y$, in some sense. However, for the purposes of this current work, we will ignore this similarity criterion and search for `simple' candidates.

\begin{figure}
  \includegraphics[width=0.49\textwidth]{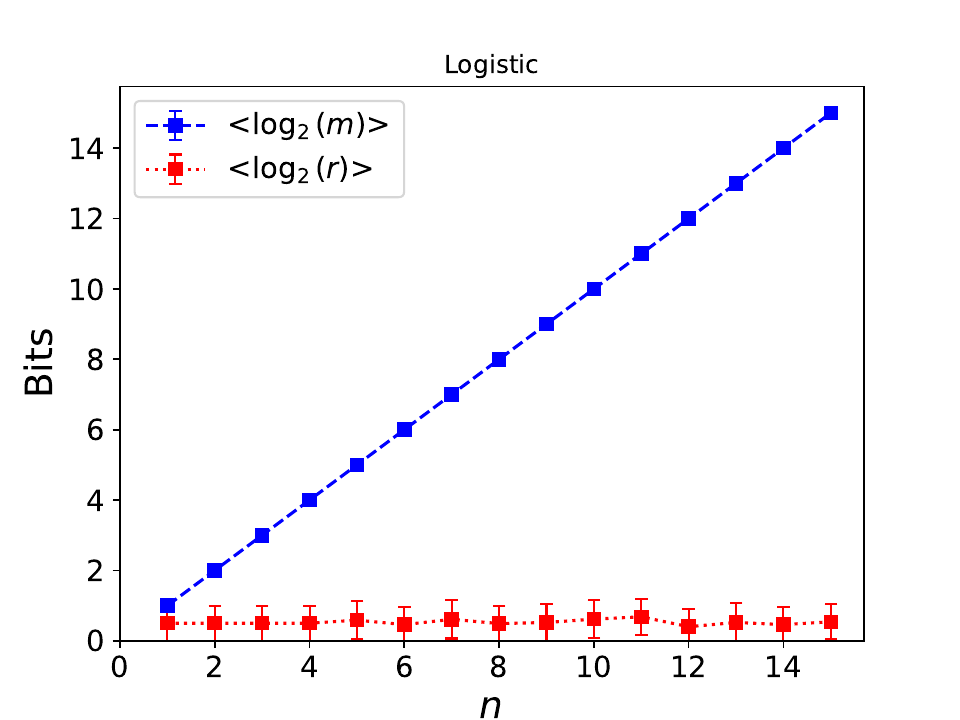}
  \includegraphics[width=0.49\textwidth]{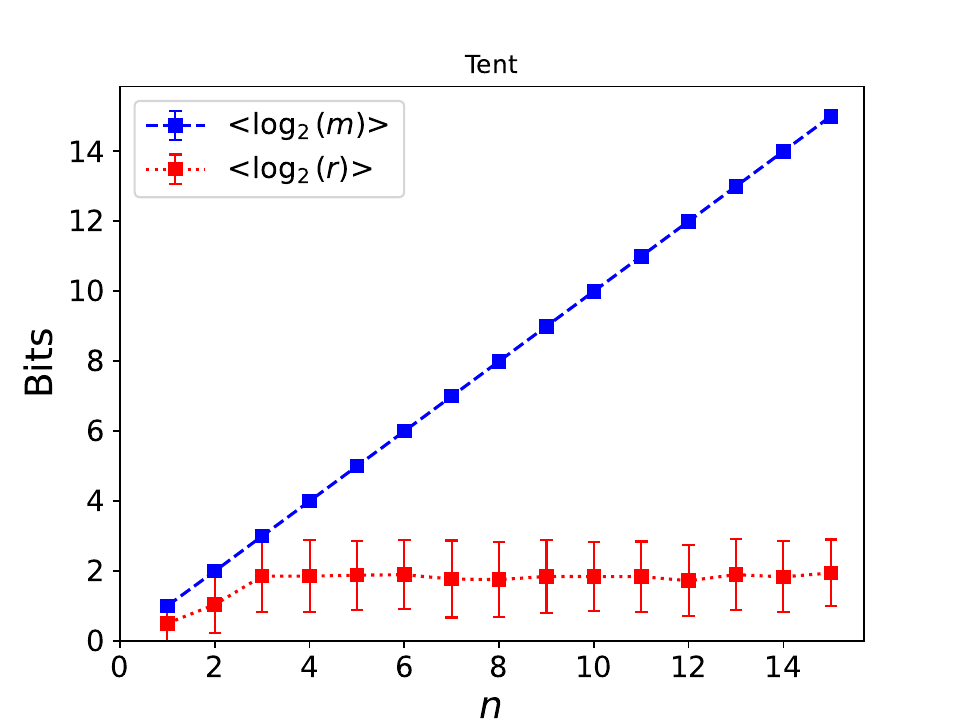} \\
  \includegraphics[width=0.49\textwidth]{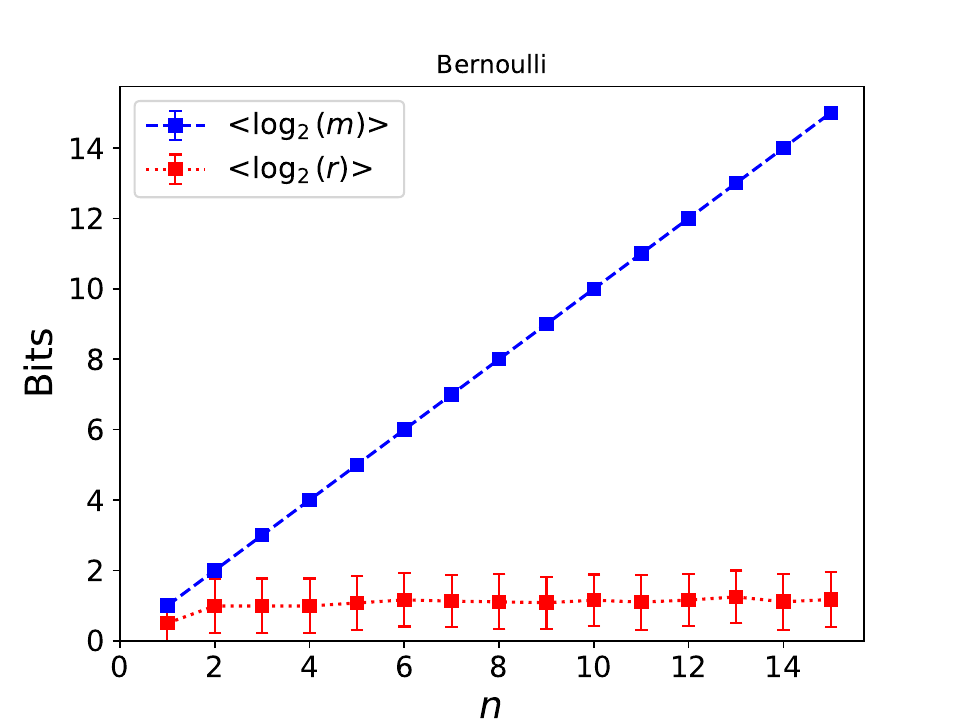}  
    \includegraphics[width=0.49\textwidth]{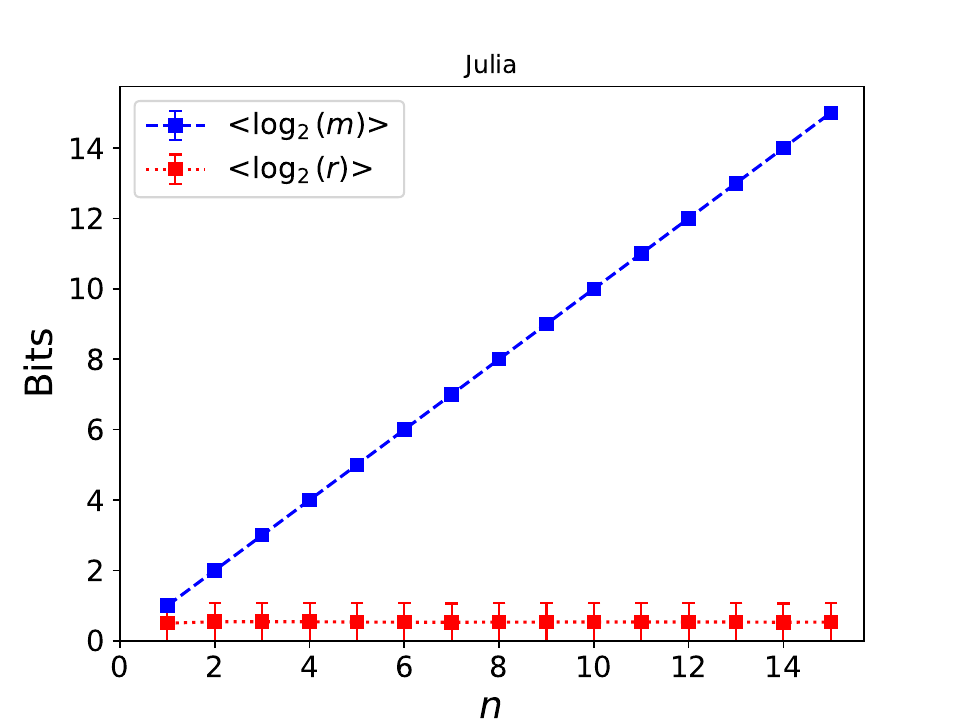} 
   
  \caption{
  Retrodiction with digit complexity method ($D\leq 2$). (\emph{top left}) logistic map; (\emph{top right}) Tent map; (\emph{bottom left}) Bernoulli map; (\emph{bottom right}) Julia map.  In each panel, the mean number of distinct candidates $m$ rises exponentially (blue), typically as $m=2^n$ (error bars $\pm$1 std.\ dev.). In contrast, the mean log rank $r$ (red) of the true starting value has a much lower value and does not continue to grow. Hence, $\log_2(r)\ll \log_2(m)$, typically, for large enough $n$, and the true starting value is in a small subset of the candidates.
    }\label{fig:digitcomplexity}
\end{figure}

\section{Experiments}

In the preceding section, we argued from two perspectives that searching for low-complexity candidates is a good way to search for the true starting value $x$. We now numerically test this candidate-complexity method.

\subsection{Applications to the maps}

\paradot{Logistic map}
The logistic map was previously studied in the context of retrodiction by Rupprecht and Vura \cite{rupprecht2018limits}.
In their work, the authors examined a distribution over possible starting values, rather than a fixed starting value that we consider. Also, while they also employ an information theoretic approach, they use Shannon information theory and entropy rather than algorithmic information theory. So while the context of their and our problem is similar, the set-up and methods are quite dissimilar.

We tested the digit complexity method described above for the logistic map, using iteration numbers $n=1, 2, 3,\dots, 15$. We used 100 starting values $x_0 $ sampled from values in [0,1] with $D=2$ (e.g., 0.45) or $1$ (e.g., 0.4) decimal places. The reason for using small $D$ at this point is to limit the complexity of the starting values (later $D$ will be increased). 
Let $r$ be the rank of the true starting value so that  $1\leq r \leq m$. (Because some complexity values may be repeated, ranking by complexity may yield ties. However, the ties can be broken by ordering candidates with equal ranks by, for example, numerical size.) For each iteration number $n$, we calculated the total number of distinct candidate values $m$, which is typically $m=2^n$, but sometimes $m$ drops due to repeated values. The mean value $<\log_2(m)>$ is taken over all starting values, and is typically very close to $n$. Additionally, for each of the different starting values $x_0$ we found the $\log_2(r)$ value, and then found the mean of these log ranks, denoted as $<\log_2(r)>$. Figure \ref{fig:digitcomplexity} shows  the  plots; as is clear $<\log_2(r)>$ is much lower than $<\log_2(m)>$ even for small $n$. Also, as expected from the above arguments, the value of $<\log_2(r)>$ does not continue to rise arbitrarily. Hence the number of bits saved is unbounded as $n$ grows, i.e., $\log_2(m) - \log_2(r)\rightarrow \infty$.

This exact retrodiction method works very well in this experiment. However, the explanation for the high accuracy is the somewhat trivial reason that the other candidates are irrational numbers with $D=\infty$, and because these irrationals are pseudo random numbers which standard lossless compression algorithms cannot handle (e.g., \texttt{zlib.compress}), then these candidates will assign high complexity to these values.

\paradot{Tent map}
The Tent map is a harder retrodiction problem as compared to the logistic map for two reasons: one is that the forward and backward iterates are uniformly distributed over the [0,1] interval, hence maximising our ignorance of starting values. The other reason is that the pre-images do not contain irrational numbers which cause $D=\infty$, but rather $D$ can be small for many candidates, and indeed even smaller than that of $x_0$. That is, the list of candidates may contain candidates which are simpler than $x_0$. To illustrate the complexity method in more detail, in Appendix \ref{app:example_tent} we show a fully worked example for the Tent map with $x_0=0.68$, $n=7$. 

We tested the digit complexity method in the same way as the logistic map, using the same values of $n$ and $x_0$. See Figure \ref{fig:digitcomplexity} which shows that $<\log_2(r)>$ is much lower than $<\log_2(m)>$ even for small $n$. Also, as expected from the above arguments, the value of $<\log_2(r)>$ does not continue to rise arbitrarily. Unlike the logistic map, in the Tent map, the red curve rises roughly linearly for $n\leq 3$, which is due to the presence of several simpler candidates which push the rank of $x_0$ higher.

The digit complexity method is successful in this example, and it is not due to the somewhat trivial reason for success in the logistic case (i.e., that $D=\infty$ for most candidates).

\begin{figure}
    \includegraphics[width=0.49\textwidth]{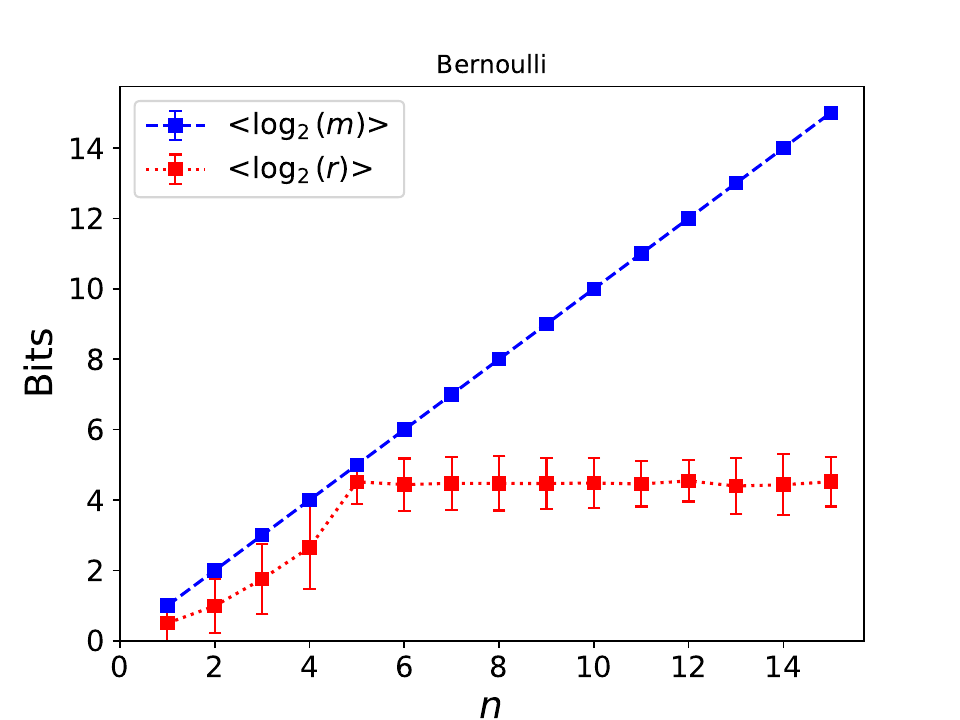}
     \includegraphics[width=0.49\textwidth]{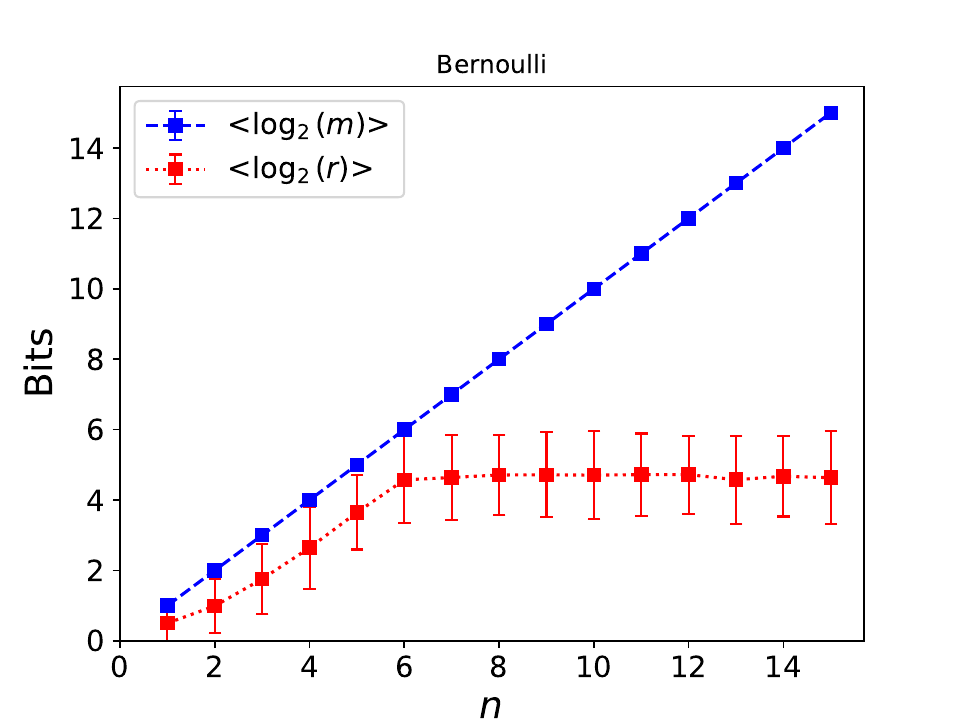}
      \includegraphics[width=0.49\textwidth]{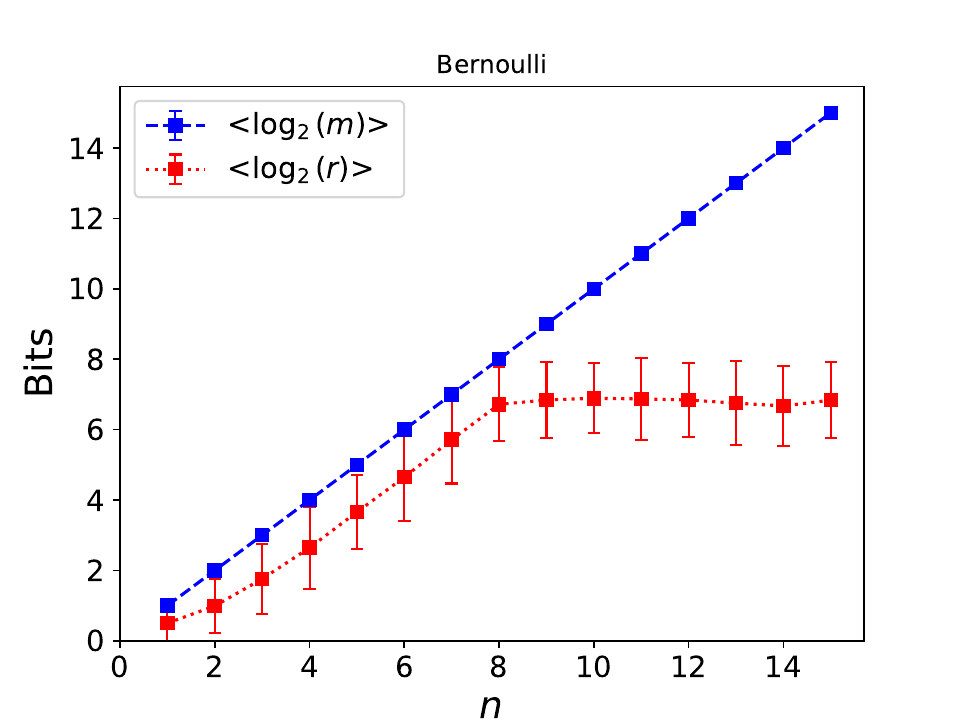}
       \includegraphics[width=0.49\textwidth]{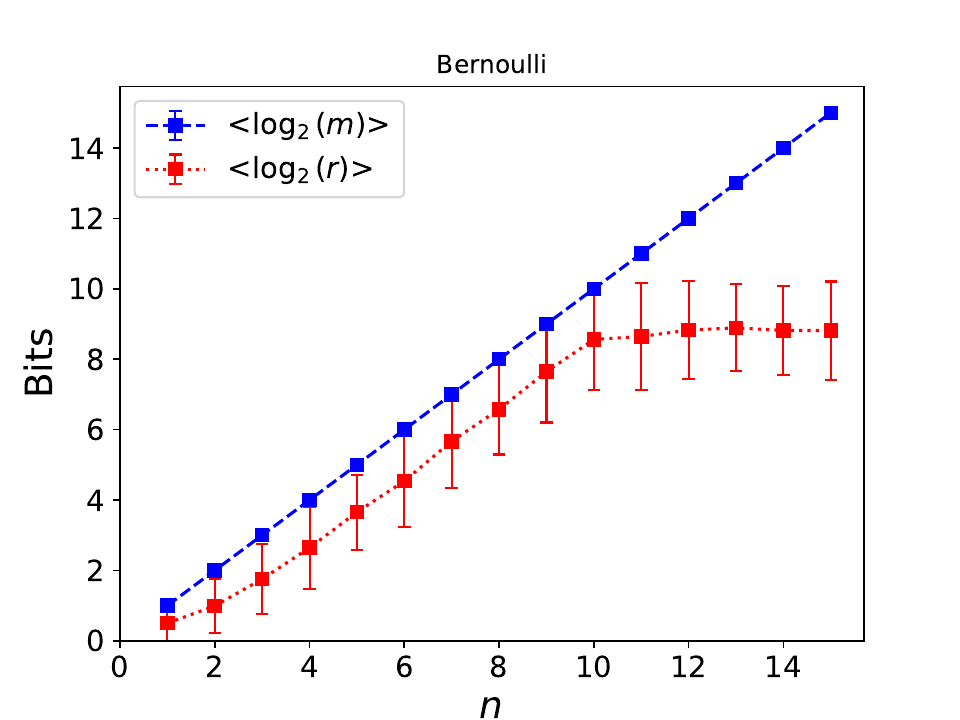}
  \caption{
  Bernoulli map retrodiction with higher decimal places, $D$. (\emph{top left}) $D=4$; (\emph{top right}) $D=6$; (\emph{bottom left}) $D=8$; (\emph{bottom right}) $D=10$.  As $D$ gets larger, typical complexity values increase, hence attaining accurate retrodiction requires larger and larger $n$ values. Each data point is an average over 50 samples ($\pm$ 1 std. dev.).
    }\label{fig:digitcomplexity_complex_x0}
\end{figure}

\paradot{Bernoulli map}
The Bernoulli map is quite similar to the Tent map, and similarly it is more challenging to do retrodiction in this map as compared to the logistic case above, and for the same reasons.

We tested the digit complexity method in the same way as the logistic map, using the same values of $n$ and $x_0$. See Figure \ref{fig:digitcomplexity} for plots of $<\log_2(m)>$ and $<\log_2(r)>$ which show that bits are saved, even for small $n$. Also, as expected from the above arguments, the value of $<\log_2(r)>$ does not continue rise. The method works well for in this experiment.

\paradot{Julia/Mandelbrot map}
Finally, looking at the Julia map, we performed numerical experiments where the starting value is $z_0=a+ib$, and $a$ and $b$ can take values -0.05, -0.025, 0.025, or 0.05, and $c=c_1 + c_2 i$ can take values -0.03, -0.01, 0.01, 0.03. This yields $4^4=256$ different retrodiction problems. In this experiment, we used very small values of $a$ and $b$ as well as small values of $n$, in order to reduce numerical issues (namely that the true starting value does not appear in the candidates).

Choosing $n = 1, 2, 3, \dots, 15$ iterations, this implies up to $m=2^{15}=32768$ candidate values for $z_0$, so that $\log_2(m)=1, 2, 3, \dots, 15$. Running the 256 retrodiction problems yields that for nearly all of them, the correct $z_0$ was ranked first in the list, $r=1$. See Figure \ref{fig:digitcomplexity}. Hence, in this map, studying the complexity of candidates via counting digits works well.

\subsection{Higher complexity starting values}

Above we saw that for very low complexity starting values such as 0.1 and 0.35 with $D=1$ and $2$, respectively, searching for simpler starting values works well. Now we consider more complex starting values. According to the theory developed above, the retrodiction method should work for both simple and complex starting values, but in the latter case it is expected that the number of iterations $n$ must be larger before $x_0$ is identifiable due to its relative simplicity.

We briefly explore how the complexity of the starting values affects the number of bits saved. We study just one map, namely the Bernoulli map, while increasing the number of digits in $x_0$ gradually, i.e., $D=4,6,8,10$. Increasing the number of digits is a proxy for increasing the complexity of starting values, because typical random starting values with larger $D$ will be more complex than with smaller $D$. 
Figure \ref{fig:digitcomplexity_complex_x0} shows that as $D$ becomes larger (i.e., $x_0$ becomes more `complex' on average), $n$ needs to be larger before saving bits. This accords with the expectations from the above theory.

\subsection{Comparison to other retrodiction methods}

In this section we describe some other retrodiction inference methods which can be used to compare to our proposed AIT-based method. 

\subsubsection*{Random.}
This method, which acts as a baseline, enumerates all the $2^n$ preimages, and then just guesses one candidate at random as the retrodiction method.

\subsubsection*{Gaussian Process (GP) beam search.}

The core idea of the GP-based method is to efficiently explore the enormous space of possible past trajectories without exhaustively evaluating each one, which would be computationally taxing. Gaussian Processes (GPs) are employed to approximate the local inverse dynamics of the system, providing probabilistic estimates of each preimage and its associated uncertainty. Because these GP predictions are only approximate, Newton iterations are applied to project each estimate exactly onto the dynamical constraint 
\begin{equation}
    f(x_t) = x_{t+1},
\end{equation}
ensuring that the reconstructed states lie on the true manifold of the map.

Each candidate history is then scored using a maximum-likelihood criterion that combines \emph{GP sharpness} (favoring low predictive variance) and \emph{Jacobian-based penalties} (discouraging trajectories in regions of strong local expansion). The resulting retrodicted trajectory is therefore both dynamically consistent and statistically most likely to be shadowed by the true path.

Importantly, the GP variance depends on the number of forward or backward iterations, $n$, since uncertainty accumulates as predictions are composed over multiple steps. This growing variance must be explicitly accounted for when ranking candidate trajectories to prevent overconfidence in long retrodictions.

Finally, a beam search is used to retain only the most promising candidate sequences at each step, allowing the algorithm to approximate the local inverse dynamics efficiently while balancing computational tractability and physical consistency. Here's a breakdown of how it works:

\paragraph{Modeling the Inverse Dynamics.}
A GP is trained to approximate the inverse mapping from the current state $x_{t+1}$ to possible predecessors $x_t$.
Instead of directly inverting the forward map $f$, which is generally infeasible for chaotic systems, the GP predicts a distribution over candidate preimages:
\begin{equation}
    x_t \sim \mathcal{N}(\mu(x_{t+1}), \sigma^2(x_{t+1})).
\end{equation}
Here, $\mu(x_{t+1})$ is the GP mean prediction (the most likely preimage), and $\sigma^2(x_{t+1})$ quantifies the predictive uncertainty.
Low-variance regions are treated as more reliable, guiding the search toward dynamically consistent regions of the state space.

\paragraph{Beam Search over Candidate Trajectories.}
Since each state generally has multiple preimages, tracing backward $n$ steps generates an exponentially large search tree.
To keep computation tractable, the algorithm performs a \emph{beam search}: at each backward step, it retains only the top $B$ most promising partial trajectories, where $B$ is the beam width.

Each candidate trajectory is scored according to two main criteria:
\begin{itemize}
    \item \textbf{GP Confidence:} candidates with low predictive variance $\sigma(x_t)$ are preferred;
    \item \textbf{Dynamic Stability:} trajectories with smaller local expansion rates (measured via $|f'(x_t)|$) are favored, as they are more likely to shadow true trajectories.
\end{itemize}
The combined trajectory score is defined as
\begin{equation}
    S = -\log \sigma(x_t) + \lambda \log |f'(x_t)|,
\end{equation}
where $\lambda$ balances the trade-off between GP certainty and dynamical stability.

\paragraph{Newton Refinement for Dynamic Consistency.}
The GP predictions $\hat{x}_t$ are only approximate preimages, typically satisfying $f(\hat{x}_t) \approx x_{t+1}$.
To enforce exact \emph{dynamical consistency}, each candidate is refined using Newton’s method, solving
\begin{equation}
    f(x_t) - x_{t+1} = 0.
\end{equation}
Starting from the GP estimate $\hat{x}_t^{(0)}$, we iteratively update
\begin{equation}
    \hat{x}_t^{(k+1)} = \hat{x}_t^{(k)} - 
    \frac{f(\hat{x}_t^{(k)}) - x_{t+1}}{f'(\hat{x}_t^{(k)})},
\end{equation}
and in higher dimensions,
\begin{equation}
    \hat{\mathbf{x}}_t^{(k+1)} = 
    \hat{\mathbf{x}}_t^{(k)} - 
    [J_f(\hat{\mathbf{x}}_t^{(k)})]^{-1} 
    \big(f(\hat{\mathbf{x}}_t^{(k)}) - \mathbf{x}_{t+1}\big),
\end{equation}
where $J_f$ is the Jacobian of $f$.
Iterations continue until
\begin{equation}
    \| f(\hat{x}_t^{(k)}) - x_{t+1} \| < \varepsilon,
\end{equation}
or a maximum number of steps (typically 3--5) is reached.
This refinement projects each approximate point back onto the manifold of exact trajectories, ensuring that the final candidates strictly satisfy $f(x_t) = x_{t+1}$.

\paragraph{Summary of the Procedure.}
At each backward step:
\begin{enumerate}
    \item The GP proposes candidate preimages for the current state;
    \item Each candidate is scored by GP uncertainty and local Jacobian magnitude;
    \item The top $B$ candidates are retained in the beam;
    \item Newton’s method refines each candidate for dynamical consistency;
    \item The process repeats until the initial state $x_0$ is reached.
\end{enumerate}
The final retrodicted trajectory is the one with the highest cumulative score, balancing statistical plausibility and exact dynamical validity.

\paragraph{Discussion.}
The GP beam search combines \emph{statistical inference} (through GP-based predictions) and \emph{deterministic correction} (via Newton refinement).
It efficiently explores the most likely backward trajectories without enumerating the full preimage tree.
However, our experiments show that even this hybrid approach performs no better than random selection in exact retrodiction tasks, underscoring the intrinsic difficulty of reconstructing past states in chaotic systems.

\subsubsection*{High-density baseline.} 
The density method ranks preimages based on the local density of candidates in their neighborhood. The intuition is that true preimages are potentially more likely to occur in regions where candidate preimages are clustered closely together. The method assigns higher scores to candidates with smaller nearest-neighbor distances, and selects the densest point. 
\begin{itemize}
    \item Rank candidates by local density (inverse nearest-neighbor distance).
    \item Approximate; requires enumeration of all $2^n$ preimages.
    \item Pros: simple, nonparametric, reveals structure in preimage clusters.
    \item Cons: exponential cost, sensitive to neighborhood choice, ignores dynamics.
    \item Use for small $n$ or as an interpretive baseline for structural clustering.
\end{itemize}

\subsubsection*{Particle filtering.} 
The particle-based method samples a large number of hypothetical initial states and propagates them forward. Each particle is weighted by how well it matches the observed $x_n$, given an assumed observation noise model. The posterior distribution over particles represents an approximate distribution of past states. The estimate is obtained from the weighted average or from the maximum a posteriori particle. 
    \begin{itemize}
      \item Samples full posterior over $x_0$ and trajectory.
      \item \emph{Principled}, handles noise explicitly.
      \item Pros: full uncertainty quantification.
      \item Cons: computationally heavy; tuning required.
      \item Use when uncertainty estimates are important.
    \end{itemize}

\subsubsection*{AIT method using different bases and complexity measures.} 
The AIT-inspired candidate-complexity method above ranked candidates by complexity while they were in decimal form, but we can also compare to how the method works when the candidates are given in binary form. 

    \begin{itemize}
      \item Choose past with lowest complexity.
      \item Pros: conceptually simple; no noise model needed.
      \item Cons: only surrogate for true Kolmogorov complexity.
    \end{itemize}

The method above used \texttt{zlib}, but there are other complexity estimators available, such as the Block Decomposition Method (BDM) \cite{zenil2018decomposition}, which is intended as an approximation to Kolmogorov complexity based on algorithmic probability. We will compare three AIT methods, namely \texttt{AITDEC} which uses decimal numbers and \texttt{zlib} as above; \texttt{AITBIT} which uses binary numbers and \texttt{zlib}; and  \texttt{AITBDM} which uses binary numbers and the BDM complexity measure. (Note that the BDM method in Python does not currently accept decimal numbers with 10 symbols for compression.)

\subsubsection*{Shadowing-inspired optimization.} 
A central concept in dynamical systems theory is the \emph{shadowing property}. 
In chaotic systems, an approximate trajectory (a \emph{pseudo-orbit}) that accumulates 
small local errors can often be ``shadowed'' by a true orbit that remains close over 
the entire time interval. More precisely, if $\{y_t\}_{t=0}^n$ is a 
$\delta$-pseudo-orbit satisfying
\[
  |y_{t+1} - f(y_t)| < \delta, \qquad t=0,\dots,n-1,
\]
then under the shadowing lemma there exists a true trajectory 
$\{x_t\}_{t=0}^n$ of the map $f$ such that 
$|x_t - y_t| < \varepsilon$ for all $t$, where $\varepsilon$ depends on $\delta$.
This property provides a theoretical guarantee that approximate backward orbits 
constructed from noisy data can, in principle, be closely followed by exact solutions.

The GP method and the AIT-based approach developed here are \emph{shadowing-inspired}. They construct approximate 
backward trajectories from a terminal observation and then refine or select among 
them in such a way that they are plausible candidates to be shadowed by a true orbit.

\begin{figure}[t]
\centering
\begin{tikzpicture}[
    >=Latex,
    every node/.style={font=\small},
    true/.style={draw=blue!60, very thick},
    pseudo/.style={draw=red!70, dashed, thick},
    dot/.style={circle, fill=black, inner sep=1pt}
]

\coordinate (x0) at (0,0);
\coordinate (x1) at (1.6,0.5);
\coordinate (x2) at (3.2,0.9);
\coordinate (x3) at (4.8,1.1);
\coordinate (x4) at (6.4,1.25);

\coordinate (y0) at (0,0.1);
\coordinate (y1) at (1.6,0.6);
\coordinate (y2) at (3.2,1.1);
\coordinate (y3) at (4.8,1.4);
\coordinate (y4) at (6.4,1.6);

\draw[true] (x0) -- (x1) -- (x2) -- (x3) -- (x4);
\foreach \p/\lab in {x0/$x_0$,x1/$x_1$,x2/$x_2$,x3/$x_3$,x4/$x_4$}
  \node[dot,label=below:\lab] at (\p) {};

\draw[pseudo] (y0) -- (y1) -- (y2) -- (y3) -- (y4);
\foreach \p/\lab in {y0/$y_0$,y1/$y_1$,y2/$y_2$,y3/$y_3$,y4/$y_4$}
  \node[dot,fill=red!70,label=above:\lab] at (\p) {};

\foreach \p in {x0,x1,x2,x3,x4}
  \draw[blue!40, thick] (\p) circle (0.25);

\draw[<->,gray] (x2) -- (y2)
  node[midway,right] {$\varepsilon$};

\end{tikzpicture}
\caption{Illustration of the shadowing property. 
A pseudo-orbit $\{y_t\}$ (red, dashed) with small step errors is 
\emph{shadowed} by a true orbit $\{x_t\}$ (blue, solid), 
which remains within distance $\varepsilon$ of it for all $t$.}
\label{fig:shadowing}
\end{figure}

\begin{figure}
    \includegraphics[width=0.9\textwidth]{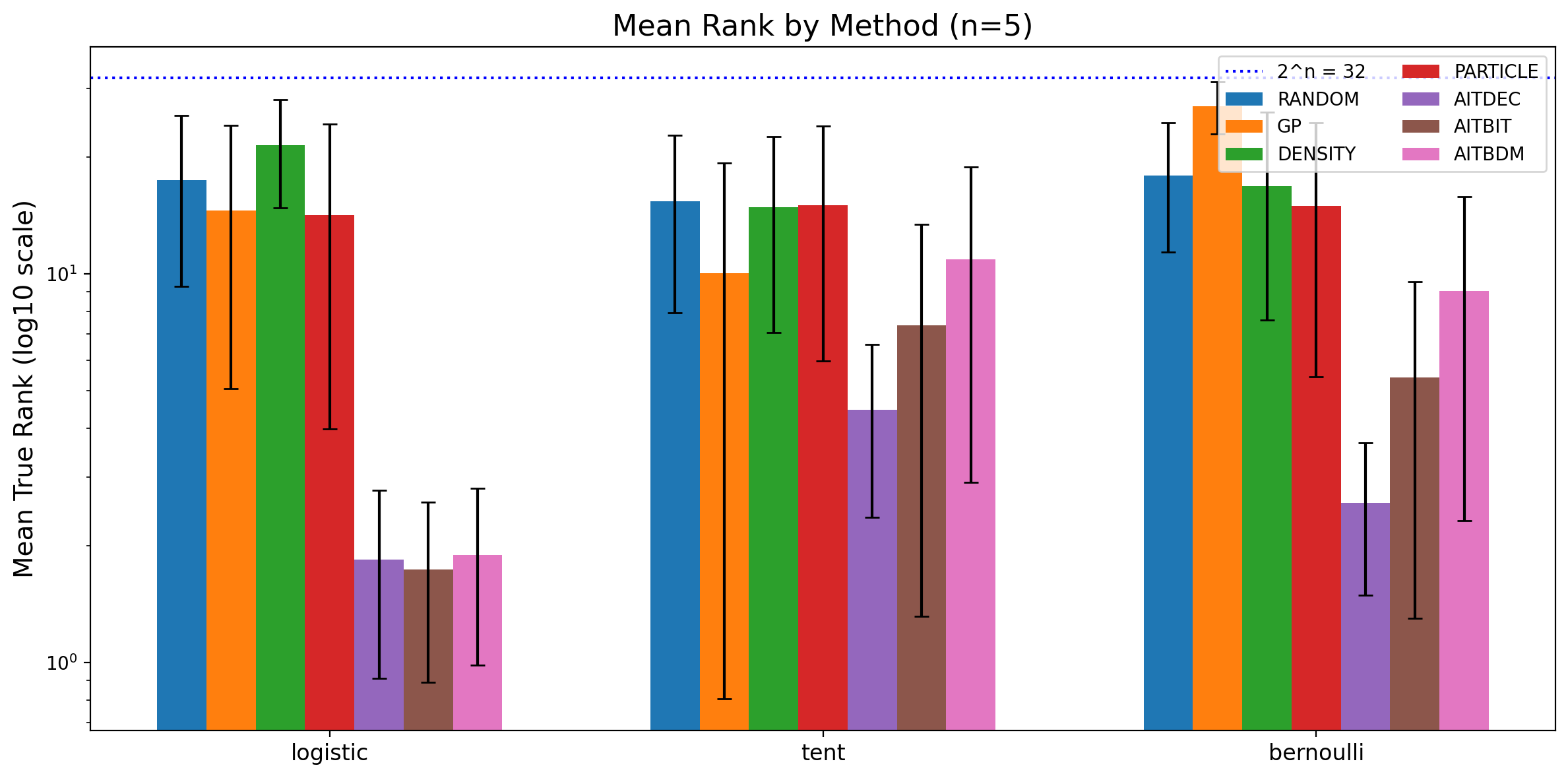}
     \includegraphics[width=0.9\textwidth]{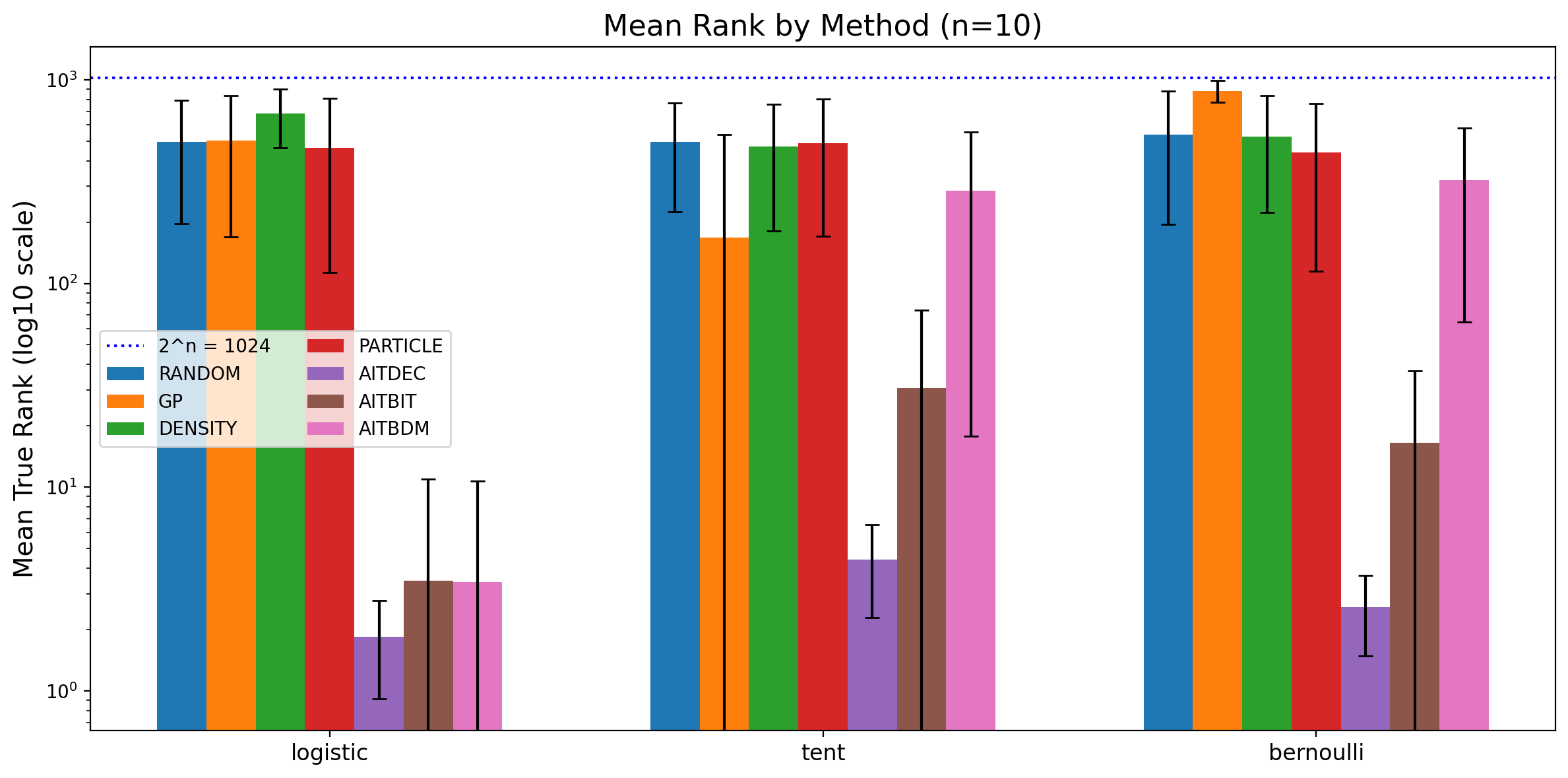}
  \caption{
   A selection of comparison retrodiction methods do no better than mere random guessing as a method of retrodiction. All three AIT approaches outperform the comparison methods. Bar heights are mean rank ($\pm$ 1 std. dev.). The top panel is for $n=5$ iterations, and bottom panel is for $n=10$ iterations.
    }\label{fig:comparison_methods}
\end{figure}

In Figure \ref{fig:comparison_methods} we compare the above methods for the three 1D maps, i.e. Tent, Bernoulli, and logistic. The vertical axes shows the mean ranks of the true starting value, as determined by each method just outlined. The blue horizontal line gives $2^n$. Methods with bars much lower than the blue dotted line are deemed more accurate. We use $n=5$ and $n=10$ iterations, and in both cases, we see all comparison methods perform no better than the random guessing. This highlights the difficulty of the problem we study, retrodiction of chaotic systems. On the other hand, the AIT-based methods, whether using decimal or binary, or the BDM complexity measure, outperform the other comparison methods. The \texttt{AITDEC} method is best for this task.

\begin{figure}
    \includegraphics[width=0.7\textwidth]{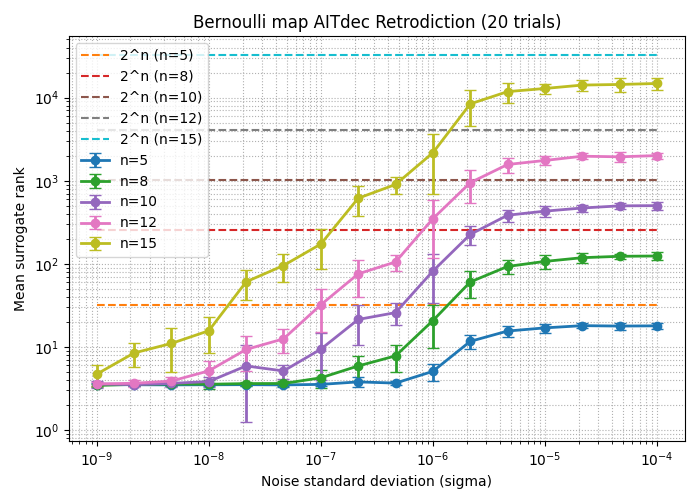}
     
  \caption{
   Effects of noisy trajectories on retrodiction accuracy for the Bernoulli map. Increasing the noise level decreases the accuracy, especially for larger $n$. For very small noise levels, and small $n$ partial retrodiction is still possible.
    }\label{fig:noisy}
\end{figure}

\subsection{Additive noise}

In the discussion above on the limitations of the AIT-based method, we noted that the map is assumed to deterministic and all candidates are known precisely.  Real-world problems in e.g. physics are typically noisy with stochastic perturbations, however, which implies a limitation on the applicability of the method so-far discussed. 

Looking more carefully, we see that for small $n$ and very small noise levels, adding random values $\epsilon_i$ to the trajectories $x_0$, $x_1$, $\dots$, i.e. $f(x_{i+1})=x_i + \epsilon_i$ will not completely disrupt the method. To see why, note that a decimal number 0.250000000000 or 0.33333333333 are both `simple' (compressible) even if very small noise is added such that only the latter digits are disrupted (e.g., they become 0.250000000007 or 0.33333333331). However, if the digits are stored to infinite precision, and the noise is also added to infinite precision, then even adding small noise will disrupt the method.

In the case where there is no noise, then $x_0$ will be in the list of candidate starting values, and hence we can investigate its rank $r$ in the complexity-sorted candidates (as above). In the case where we add noise, $x_0$ will not be in the list of candidates. (Note that the list of candidates is deterministically inferred from the noisy $x_n$.) However, if the noise is small enough, there will be an element, $x_0'$, which is very close to $x_0$, and $|x_0-x_0'|\approx 0$. Assuming this, we can investigate the rank of $x_0'$ in the list of noisy candidates; we will call this the surrogate rank. If the surrogate rank is low, $r\approx 1$, and $|x_0-x_0'|\approx 0$ is small, then this will count as successful retrodiction, because we can identify $x_0'$ with non-trivial accuracy, and $x_0'$ is a good estimate of $x_0$.

In Figure \ref{fig:noisy} we perform a brief experiment of the effect of noise on retrodicting the Bernoulli map. We use additive noise given by a Gaussian distribution with standard deviation $\sigma$ ranging from $10^{-9}$ to $10^{-4}$, which is added in each forward iterate. We experiment with $n=5, 8, 10, 12, 15$. For noise level of $10^{-9}$, successful retrodiction is achieved for all $n$ values, as evidenced by observing that the (surrogate) rank is much lower than $2^n$. Increasing the noise to level $10^{-8}$ raises the rank (decreases accuracy) while for $n\leq 12$ the accuracy is maintained. Increasing the noise level more decreases the accuracy, especially for larger $n$, as is intuitive.

\section{Approach 2: Nearest-neighbor-based reconstruction }\label{sec:nnrecon}

\subsection{Approach details}

The complexity of candidates method above has been shown to generally effective, but at the same time it has some important limitations for practical retrodiction problems (detailed in Section \ref{limitations}). Here we propose another method which is more practical, and yet still motivated by AIT arguments.
The following could be regarded as a more abstract version of the complexity of digits method: 
For each $x^i\in\cX_0\subseteq\mathbb{R}^d$ (not to be confused with the $i$th iterate $x_i$) 
find the distance $S^i$ to its nearest neighbor in $\cX_0\setminus\{x^i\}$.
This distance $S^i$ could be used to estimate the ``denseness'' 
of $\cX_0$ at $x^i$ as $\rho(x^i)\propto (S^r)^{-d}$%

From another perspective, `simpler' candidates are supposed to be different from the remaining typical members (in some sense), and so we might expect these `simpler' candidates to be outliers (see below for more on this), e.g., in low-density areas. One advantage of a density-based method is that it does not require knowing the details of the map $f$. Indeed, directly from the dataset, one can search for low-density regions, specifically the candidates in those regions.

To proceed, w.l.g.\ assume $x=x^1$. If $S^1$ is large, we only need to encode $x$ very crudely as $\hat x$, 
since we can reconstruct $\cX_0$ exactly from $y$ and $g$ and then search for the closest point to $\hat x$ in $\cX_0$. 
Roughly, we need to encode each dimension of $x$ to $\log_2(1/S^1)$ bits.
More precisely, we convert $S^r$ to a probability distribution 
$p_{x^r}:=\fr1c(S^r)^d$, where $c:=\sum_{r=1}^m (S^r)^d$ which leads to code length
$\log_2 p_{x^1}^{-1}=d\log_2(1/S^1)+\log_2 c$.

Alternatively we can sort $\cX_0$ in decreasing order of $S^r$, 
and look for the rank $r$ of $x$, and then encode the rank in CL$(r)$ bits.
This is reminiscent of the Loss Rank principle \cite{Hutter:07lorp,Hutter:10lorpx}.
This will lead to shorter code than naive $\log_2(m)$ bits if we expect $x_0$ to be in a low-density region.
For chaotic maps and especially for large $d$ this is plausible. While both the probability distribution and rank options are possible, in this current study we will use the second option of ranking the candidates (as we did above for the complexity-based method).

\begin{figure}
    \includegraphics[width=0.49\textwidth]{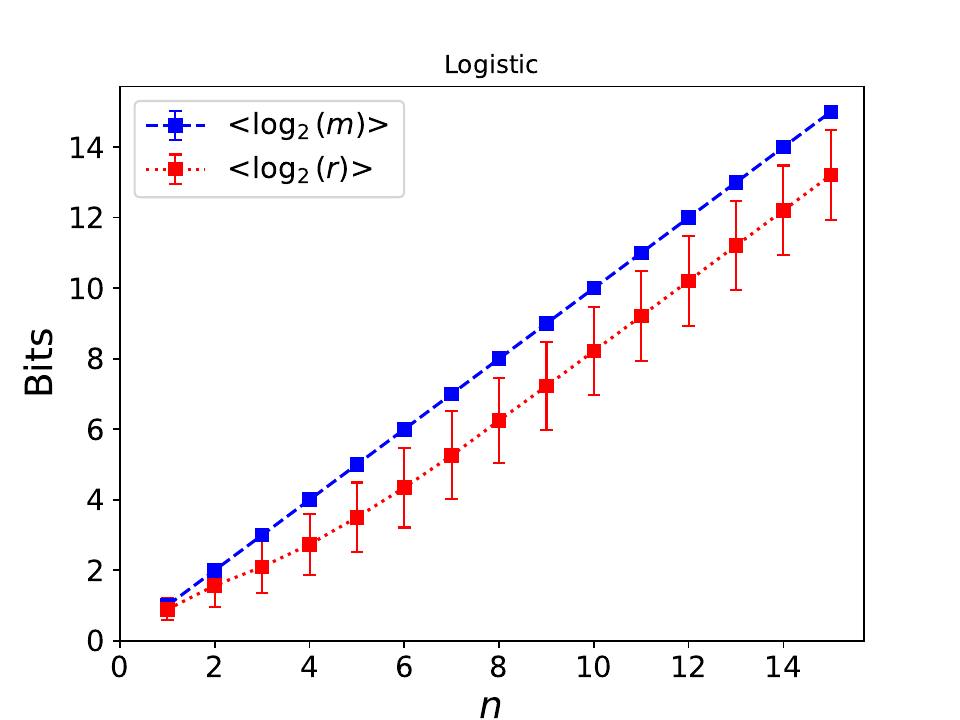}
     \includegraphics[width=0.49\textwidth]{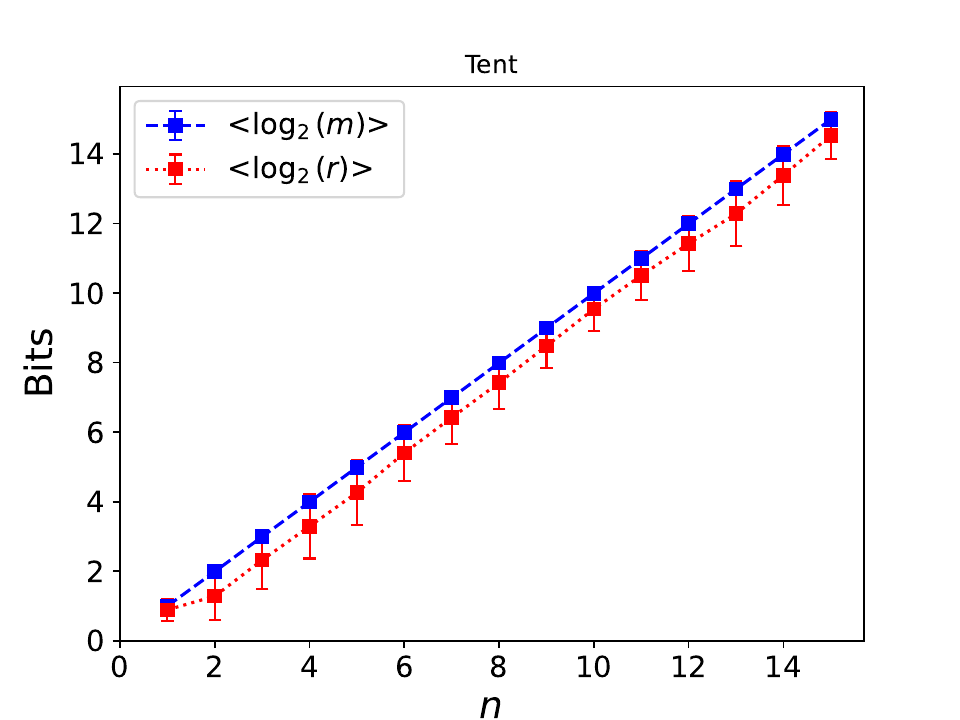}
     \includegraphics[width=0.49\textwidth]{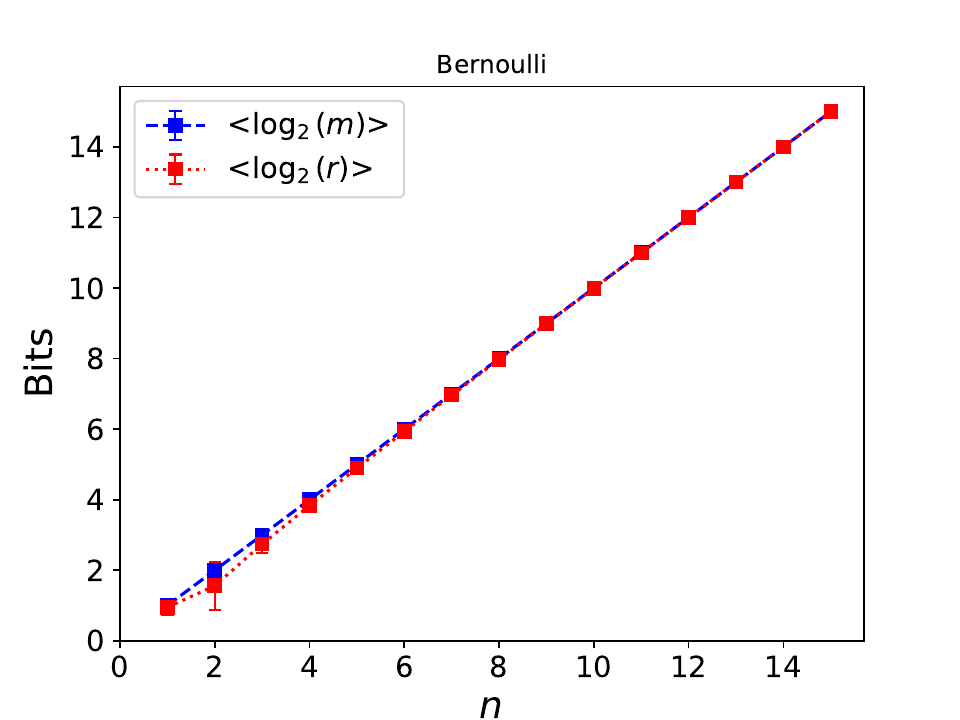}
     \includegraphics[width=0.49\textwidth]{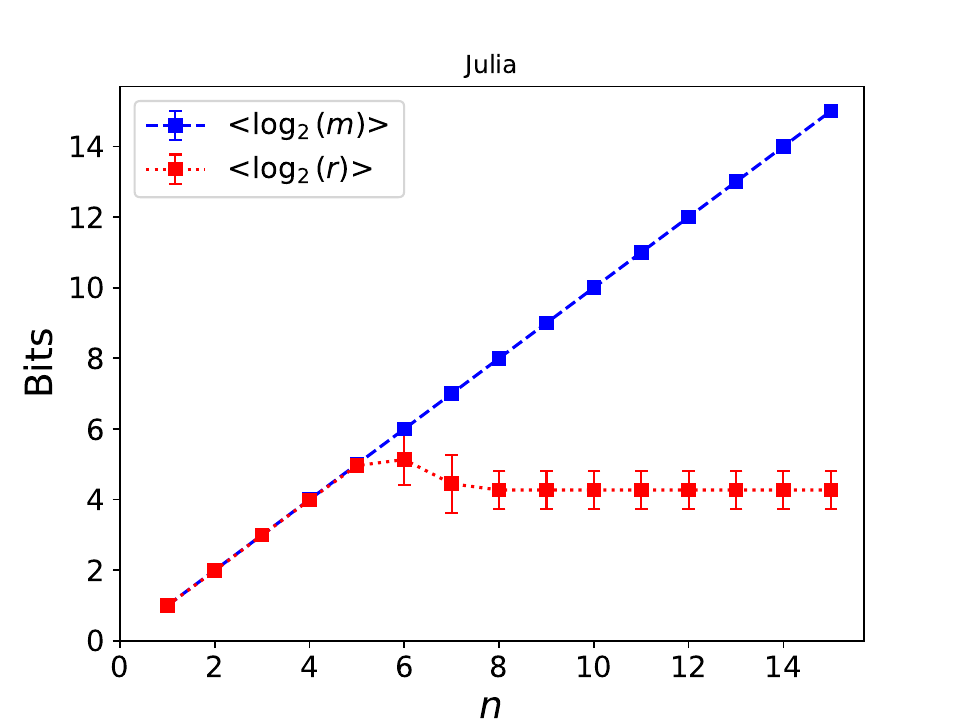}
  \caption{
  Nearest-neighbour density method. Identifying $x_0$ according to low density regions, coded by rank. For the logistic map (\emph{top left}), only a few bits are saved. For the Tent map (\emph{top right}) and Bernoulli map (\emph{bottom left}), no bits are saved (as expected). For the Julia map (\emph{bottom right}), the method works well (bits are saved) for large enough $n$.
    }\label{fig:densitymethod}
\end{figure}

\subsection{Kolmogorov complexity and outlier analysis}

Before numerically testing this second approach to the retrodiction problem, we examine the connection between low density regions and Kolmogorov complexity further. This aids in highlighting the connection of the proposed density method to the overall theme of the paper, i.e., using algorithmic information theory for retrodiction.

Searching for elements in low density regions is like searching for outliers in a data set. Moreover, there is a natural connection between Kolmogorov complexity, and outlier analysis. In the former, discrete objects are assigned complexity values, with low complexity values assigned to `simple', or rather `special' and identifiable objects. By contrast, high complexity values are assigned to random, irregular and unremarkable objects. Also, there is a  close connection between \emph{random} patterns and \emph{typical} patterns \cite{landsman2023typical}. Putting these together, we can see an intuitive connection of outliers (which are `special' data points), and inliers (which are typical data points). 

Formally, let $h$ be a function that ranks elements of a finite set by their outlier degree, with lower ranks $h\approx 1$ indicating a stronger  outlier. We will assume that the ranks are all distinct, which is not unreasonable in our setting because in real-world data, distances between data points are typically given as real numbers and hence are unlikely to match exactly, leading to non-unique rank values. Given a finite set $\mathcal{S}$ with elements $s^0,s^1,s^2,\dots, s^{N-1}$, where $N=|\mathcal{S}|$, we can form a bound connecting the Kolmogorov complexity of element $s^i$ to its outlier degree:
\begin{equation}
    K(s^i) \leq K(h(s^i)) + K(S) +O(1) 
\end{equation}
which follows easily by observing that we can describe any element $s^i$ by first generating the set $S$, and then specifying the rank $h(s^i)$ of the element $s^i$ (again assuming distinct ranks).
Further, for typical elements of $\mathcal{S}$, and assuming that $\mathcal{S}$ is a simple set, and $h$ is an $O(1)$ complexity function, we can approximate that
\begin{equation}
   K(s^i) \lesssim \log_2(h(s^i))
\end{equation}
Therefore, if $h(s^i)\approx 1$, $K(s^i)$ must be small. In other words, outliers must have low complexity. The converse is not necessarily true: there may be low-complexity elements which are not deemed to be outliers. Nonetheless, the upper bound theoretically supports connecting the two notions and approximating low complexity via low density (as above). 

Note that there are very many outlier analysis methods available \cite{aggarwal2017outlier}, including multivariate and sophisticated methods \cite{aggarwal2001outlier,zimek2012survey}. The nearest-neighbor approach used in this work is only one such outlier method. Any of these many outlier methods could be used in principle for the function $h$, but we do not experiment with these here because our work is not intended as an analysis of outlier detection methods.

In the case where $\mathcal{S}$ is not a simple set, for example, if $\mathcal{S}$  is made up of random numbers, then the connection of outliers and complexity is not so straightforward. We will leave this for future work because it is not relevant to the deterministic retrodiction problem set-up we analyze here.

\subsection{Experiments}

Now we test the nearest-neighbor density method just proposed, using the same experimental set-up, i.e., the same maps, sample parameter values, and same iteration values $n$.

\paradot{Logistic map}
Testing the density method on this map, we see in Figure  \ref{fig:densitymethod} that it is only modestly effective, saving only a few bits. Moreover, the value $<\log_2(r)>$ appears to increase nearly linearly for all tested values of $n$, rather than leveling off.

\begin{figure}
  \includegraphics[width=0.50\textwidth]{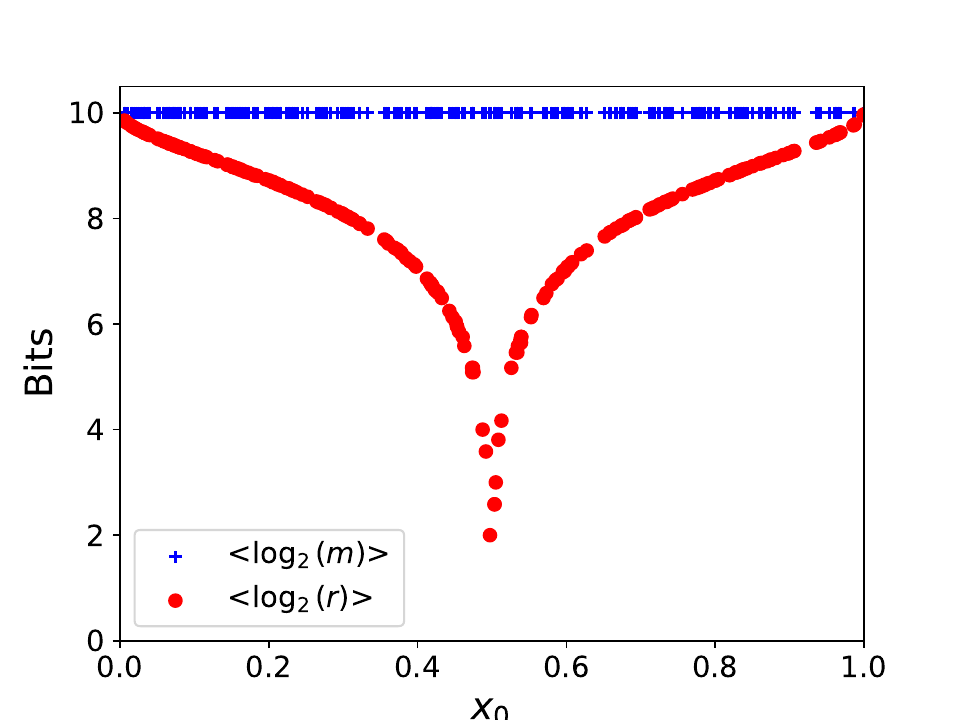}  
  \includegraphics[width=0.42\textwidth]{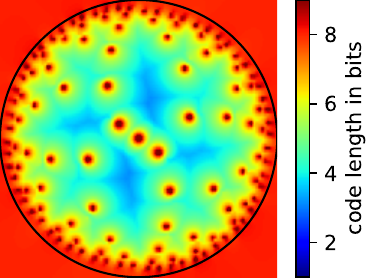}
  \caption{
    (\emph{left}) Logistic map density method.
    The code lengths for 250 random samples ($n=10$) of $x_0\in[0,1]$ together with their mean $\log_2 m$ values, and  mean $\log_2 r$ values. For $x_0\approx 0$ or 1, not many bits are saved. For $x_0\approx 0.5$, many bits are saved.
    (\emph{right})    Julia map density method. The code lengths for $z_0\in[-1,1]^2$ using $n=8$ and $c$=0.01+0.03$i$. The variation in the number of bits saved is complicated and intricate. The mean code length for $z_0\in[-1,1]^2$ is 6.21, and the mean for $|z_0|\leq 1$ is 5.66.
    }\label{fig:code_length_varying}
\end{figure}

The fact that some bits are saved at all is presumably mainly due to the fact that the testing set of starting values $x_0$ is sampled uniformly from [0,1]. This means that many values are sampled from $\approx 0.5$, values  which lie in low-density regions of the  Beta(0.5,0.5) distribution (which describes the distribution of candidates). Figure \ref{fig:code_length_varying} shows how the code lengths vary depending on the starting value, and how for $x_0\approx 0.5$ many bits are saved, evidenced by $\log_2 r \ll \log_2 m$. By contrast,  for $x_0\approx 0$ or 1, effectively no bits are saved. As expected,  more bits are saved for $x_0\approx 0.5$ where the Beta(0.5, 0.5) has low density.

\paradot{Tent and Bernoulli maps}
The Tent and Bernoulli maps repeatedly and linearly stretch out 
an $\eps$-neighborhood of any $x$ by a constant factor of $2$ independent of $x$,
so any (non-atomic) distribution converges to a uniform distribution 
in the forward and backward directions, so this density method is ineffective. In other words, all regions of the set of candidates will be roughly equally dense, so searching for the true candidate in low-density regions is not applicable.
Therefore, we do not expect this method to be effective for the Tent and Bernoulli maps. The numerical simulations illustrate: see Figure \ref{fig:densitymethod}, which shows that the density method saves no bits for these two maps, regardless of the value of $n$.

\paradot{Julia map}
The Julia has dimension $d=2$, because $z\in\mathbb{C}\equiv\mathbb{R}^2$. For this map with fixed $c\in\mathbb{C}$, 
the density method works well, provided $n$ is large enough:  see Figure \ref{fig:densitymethod}.
This map shows some curious behavior. With $n\leq 5$ no bits are saved, but for $n>5$ many bits are saved. The explanation is apparent from Figure \ref{fig:julia_n58} in the Appendix, which shows for $n=5$, the candidate points are spread roughly uniformly, such that searching for low-density regions is not helpful. However, the same figures illustrate how when $n= 7$, the true starting value $z_0$ is in a low-density region, while the majority of candidates lie in dense regions of typical values.

While the density method works well for this map on average, similarly to the logistic map, the number of bits saved varies substantially depending on the location of the starting value. As one example, merely to illustrate this, Figure \ref{fig:code_length_varying} shows results for the Julia map using $n=8$ and $c$ = 0.01+0.03$i$, for 100 values taken from  $z_0\in[-1,1]^2$, which includes the unit circle $|z_0|\leq 1$. There is a complex pattern of how many bits are saved, with small patches of values for which very few bits are saved and other regions for which many bits are saved. Also, outside of the unit circle (drawn in black), effectively no bits are saved. We conclude the density method works very well in general for the Julia map but fails for a small subset of values.

\subsection{Stretch-based reconstruction}

Finally, we briefly consider another approach to the retrodiction problem:
A related approach to density asks what is the $\eps$-smeared probability $\rho(x)$ that $y$ maps to an $\eps^2$-neighborhood of $x$. 
 The answer for $\eps\to 0$ is $\rho(x)\propto 1/g'(x)$ for every $x\in\cX_0$, where $g'$ can be computed from $f'$ via the chain rule.
Keeping $\eps$ small but finite gives a more smoothed and possibly more interpretable $\rho$.
Equivalently, in the forward direction, we could take a narrow Gaussian $x\sim\cN(x,\eps^2)$ 
around potential initial values $x$, which implies approximately $y\sim\cN(g(x),[g'(x)\eps]^2)$.
We call $g'(x)$ the \emph{stretch factor} of $x$. The idea is to 
choose $x\in\cX_0$ such that $y$ has the minimal variance = smallest $g'(x)$ = largest ``density'' $\rho(x):=1/g'(x)$.

This stretch method is related to the density method. However, it is not fundamentally a complexity-based approach.
This method would also not work for the Tent and Bernoulli maps, and requires knowing the function $g(\cdot)$ and its derivative, which is a stronger requirement than the density approach (which is relatively straightforward and does not require this additional information). Due to the similarity of the methods, the same limitations but extra requirements, and the weak connection to complexity, we will not pursue this method here. See the Appendix for more details.

\section{Discussion}\label{sec:disc}

\subsection{Main results}

We have studied the problem of exact retrodiction in deterministic chaotic dynamical systems, a problem which consists of inferring the starting value $x_0$ of series $x_0,x_1,x_2,\dots,x_n$, given only $x_n$ and the relevant function $f$. Because the number of possible historical candidates $x_0$ typically grows exponentially with $n$, and due to mixing, the problem of identifying the starting value is very challenging and indeed impossible to solve in full generality. Despite this, we have shown that by viewing the problem from the perspective of algorithmic information theory, initial progress can be made with this problem in several settings. Indeed, a main motivation of this work is to explore algorithmic information theory applications in dynamical systems \cite{wolpert2001computational,mihailovic2015novel,zenil2019algorithmic,dingle2023note,dingle2024exploring,hamzi2024simplicity}.

Our main results are firstly that provided the starting value $x_0$ is of finite Kolmogorov complexity, for large enough $n$, eventually $x_0$ will be within a set of vanishing size relative to the full set of candidates. It follows that if good approximations to the complexities are available, then $x_0$ will be in this set of low complexity, and hence can be identified more easily as compared to a uniform random guess from all possible candidates. While we measured complexity via compressing decimal and binary digits, other more sophisticated complexity estimates could be developed/used in the future, while still guided by the same method. 
Secondly, in some maps looking for $x_0$ among low-density regions of the candidate space is an effective approach, which is also motivated by Kolmogorov complexity arguments. This second approach has the significant benefit of being easy to implement and does not require knowing the details of the map $f$. 

Interestingly, despite the mixing property of chaotic dynamical systems which is usually understood as a mechanism for `forgetting' earlier values of a series, here we have shown that for finite complexity starting values, the system may not `forget' so completely.
Also interestingly, in the first approach we have not invoked Occam's razor to solve this retrodiction problem, but at the same time have arrived at the same conclusion Occam's razor would have suggested: from among the candidate values, prefer the simpler ones.

While our approach makes significant theoretical progress, it is a first step. The proposed methods are computationally expensive, requiring the full enumeration of an exponentially large candidate set, and are sensitive to noise. These limitations highlight important directions for future work, which include improving computational efficiency, extending the analysis to random dynamical systems, and further exploring the promising connections between AIT, kernel methods, and machine learning in the context of dynamical systems.

\subsection{Reframing retrodiction methods as Gaussian Process models with different priors}

In this work, we have explored the problem of retrodiction, or inferring past states from observed future states, in chaotic dynamical systems. We formulated several methods for this task, each with varying degrees of complexity and assumptions. In this section, we further develop the discussion by reframing these methods within the context of \emph{Gaussian Process (GP)-based retrodiction} with different priors, offering a unifying perspective. These methods all fall within the broader context of \emph{complexity-based priors for GPs}, which is a central theme that will be further developed in a future paper \cite{Hamzi2024}, introducing Solomonoff Gaussian Processes (SGPs).

All of these methods share the shadowing property, which serves as a unifying feature. Each approach seeks to approximate a trajectory that dynamically shadows the true trajectory of the system, meaning that they all aim to find paths that are consistent with the system's dynamical constraints.

\subsubsection*{Random Method $\rightarrow$ GP with a Uniform Prior}
The simplest baseline method involves selecting a preimage at random from all $2^n$ candidate trajectories. This can be interpreted as a GP-based retrodiction model where the prior is {uniform} across all possible preimages. In the GP framework, this translates to a {non-informative prior}, effectively assuming that every candidate is equally likely, with no knowledge of the system’s dynamics. This serves as a lower bound for the performance of more sophisticated models.

\subsubsection*{Gaussian Process Beam Search $\rightarrow$ GP with a Dynamical Prior}
In the GP-based beam search method, a Gaussian Process is used to model the inverse dynamics, and Newton's method is applied to enforce dynamical consistency. This method can be viewed as a {GP with a dynamical prior}, where the prior reflects the underlying dynamical system's evolution, such as $f(x_t) = x_{t+1}$. The GP learns the system's inverse dynamics, allowing it to efficiently explore candidate trajectories that are dynamically consistent. The use of {beam search} ensures that only the most promising candidates are evaluated, thus improving computational efficiency.

\subsubsection*{High-Density Baseline $\rightarrow$ GP with a Density-Based Prior}
The {high-density baseline} ranks preimages based on their local density, assuming that the true preimages are more likely to lie in dense regions of the candidate set. This can be formulated as a {GP with a density-based prior}, where the GP not only learns the inverse dynamics but also incorporates a {density model} that prioritizes candidates in regions of high density. The resulting GP model is therefore influenced by both the system's dynamics and the local density of the preimages.

\subsubsection*{Particle Filtering $\rightarrow$ GP with a Likelihood Maximization Prior}
In particle filtering, a set of particles is propagated backward, and each particle is weighted based on its likelihood of producing the observed data. A GP-based version of this approach would have a {likelihood maximization prior}, where the GP posterior is updated by iterating over the candidate preimages, adjusting their weights according to the likelihood of observing the final state $x_N$. This allows the GP to focus on the most likely preimages while maintaining consistency with the observed dynamics.

\subsubsection*{AIT Methods (AITDEC, AITBIT) $\rightarrow$ GP with a Complexity-Based Prior}
The {AIT-based methods} (such as AITDEC and AITBIT) rank preimages based on their algorithmic complexity. In the GP framework, these methods can be interpreted as having a {complexity-based prior}, where the GP is biased toward simpler, less complex trajectories. The \textbf{AITDEC} method, which uses decimal representations, favors candidates that are easier to represent in decimal form, while the \textbf{AITBIT} method, based on binary representations, favors candidates that are more compressible in binary form. By incorporating these complexity penalties, the GP can prioritize simpler, more probable preimages.

This formulation of retrodiction methods as \textbf{GP-based models with different priors} provides a unified perspective, where each method represents a specific choice of prior that influences the model’s assumptions about the dynamics of the system. By varying the prior, we can tailor the GP model to different types of dynamical systems, from simple random processes to highly structured, chaotic systems.

All of these methods fall within a broader framework of \textbf{complexity-based priors for Gaussian Processes}, as explored in \cite{Hamzi2024, Hamzi2024b, Hamzi2025, HamziHutterOwhadi2025_SGP}, which bridge the gap between machine learning and algorithmic information theory. In this framework, Solomonoff Gaussian Processes (SGPs) are introduced as GPs with universal priors inspired by algorithmic information theory, offering a powerful tool for modeling complex systems.

\subsection{Extensions}
This work is only the first step in addressing the retrodiction problem for dynamical systems, and we do not claim to have fully solved this challenging problem. Some future directions include investigating the logistic map for $\mu< 4$; investigating other maps, especially those in higher dimensions; using other outlier detection methods (instead of density); studying retrodiction in random dynamical systems, and finding ways to improve the computational efficiency. Because the current work is closely related to inverse problems, there may also be avenues to pursue in that direction.   

Finally, one could explore connections with kernel methods especially given their connections with AIT as established in \cite{bh_huter_kfs_ait, BridgingAIT_part2}. Kernel-based methods  \cite{CuckerandSmale} have provided strong mathematical foundations for analyzing dynamical systems \cite{lengyel2024kernelsumsquaresdata, lee2024gaussianprocessessimplifydifferential, bhcm11,bhcm1,lyap_bh,bhks,hamzi2019kernel, bh2020b,klus2020data,ALEXANDER2020132520,bh12,bouvrie2017kernel,hb17, mmd_kernels_bh, akian2022learning, hamzimaulikowhadi, 5706920,   yk1, bh_sparse_kfs, BHPhysicaD, hamzipaillet, YANG2024134192, BHPhysicaD, lee2021learning,hamzimaulikowhadi,bhkfjpl,bhkfsdes, Lee2024KernelMethodsLyapunov, lee2024kernelmethodsapproximationeigenfunctions}. It is therefore natural to consider a kernel-based version of our current work.  In particular, the integration of Machine Learning (ML) techniques, particularly those involving kernel methods rooted in AIT as illustrated in \cite{bh_huter_kfs_ait, BridgingAIT_part2}, offers a compelling avenue for future exploration. The synergy between AIT and ML can enhance the theoretical robustness and practical application of retrodiction methods in several ways:
\begin{itemize}
    \item Kernel Methods for Retrodiction: Utilizing kernel-based methods, such as Sparse Kernel Flows, could allow for a more nuanced comparison of historical candidate states by effectively capturing the underlying dynamics of the system. These methods could help in identifying the kernel that best approximates the true starting state from complex, high-dimensional data sets by choosing kernels with lowest Kolmogorov complexity.

\item Clustering and Density Estimation: Applying kernelized clustering and density estimation methods to retrodiction problems could provide a structured way to analyze the density and distribution of candidate starting values. This approach would leverage the computational efficiency of ML to manage and interpret the vast number of potential solutions, potentially offering a more scalable solution to retrodiction in complex dynamical systems.
\item Exploiting Kolmogorov Complexity: The notion of using Kolmogorov Complexity within ML frameworks, such as through complexity-based kernels introduced in \cite{BridgingAIT_part2}, could redefine how similarities between states are quantified, moving beyond traditional metrics and embedding a deeper understanding of informational content directly into the retrodiction process.
\end{itemize}

\paradot{Acknowledgments} We would like to thank Stefano Luzzatto (ICTP) for useful discussions.

\paradot{Data and code availability} There was no data associated with this study. The code generating the figures is available from \texttt{https://doi.org/10.5281/zenodo.15782271}. 

\bibliographystyle{plain}
\bibliography{RetrodictionRefs,kernels_ait_ml} 

\newpage

\appendix

\renewcommand\thefigure{\thesection.\arabic{figure}}

\section{Example of complexity method}\label{app:example_tent}
Below an example where the initial value of $x_0=0.68$ is iterated using the Tent map with $n=7$. This yields $m=128$ different possible candidate values. By sorting candidates by complexity (\texttt{zlib.compress}), the rank of $x_0$ is found to be 7 out of 128. Note that the absolute complexity values are not important, but rather their relative complexities which gives the ranks.

\begin{verbatim}

Example for Tent map 
x0= 0.68
n= 7

List of all candidates= [0.0075   0.9925   0.4925   0.5075   0.2425   
0.7575   0.2575   0.7425
 0.1175   0.8825   0.3825   0.6175   0.1325   0.8675   0.3675   0.6325
 0.055    0.945    0.445    0.555    0.195    0.805    0.305    0.695
 0.07     0.93     0.43     0.57     0.18     0.82     0.32     0.68
 0.02375  0.97625  0.47625  0.52375  0.22625  0.77375  0.27375  0.72625
 0.10125  0.89875  0.39875  0.60125  0.14875  0.85125  0.35125  0.64875
 0.03875  0.96125  0.46125  0.53875  0.21125  0.78875  0.28875  0.71125
 0.08625  0.91375  0.41375  0.58625  0.16375  0.83625  0.33625  0.66375
 0.008125 0.991875 0.491875 0.508125 0.241875 0.758125 0.258125 0.741875
 0.116875 0.883125 0.383125 0.616875 0.133125 0.866875 0.366875 0.633125
 0.054375 0.945625 0.445625 0.554375 0.195625 0.804375 0.304375 0.695625
 0.070625 0.929375 0.429375 0.570625 0.179375 0.820625 0.320625 0.679375
 0.023125 0.976875 0.476875 0.523125 0.226875 0.773125 0.273125 0.726875
 0.101875 0.898125 0.398125 0.601875 0.148125 0.851875 0.351875 0.648125
 0.039375 0.960625 0.460625 0.539375 0.210625 0.789375 0.289375 0.710625
 0.085625 0.914375 0.414375 0.585625 0.164375 0.835625 0.335625 0.664375]

List of candidate complexity values= [46, 46, 46, 46, 46, 46, 46, 46, 46,
46, 46, 46, 46, 46, 46, 46, 45, 45, 45, 45, 45, 45, 45, 45, 44, 44, 44, 
44, 44, 44, 44, 44, 45, 45, 45, 45, 45, 45, 45, 45, 45, 45, 45, 45, 45, 
45, 45, 45, 45, 45, 45, 45, 45, 45, 45, 45, 45, 45, 45, 45, 45, 45, 45, 
45, 48, 48, 48, 48, 48, 48, 48, 48, 48, 48, 48, 48, 48, 48, 48, 48, 48, 
48, 48, 48, 48, 48, 48, 48, 48, 48, 48, 48, 48, 48, 48, 48, 48, 48, 48,
48, 48, 48, 48, 48, 48, 48, 48, 48, 48, 48, 48, 48, 48, 48, 48, 48, 48, 
48, 48, 48, 48, 48, 48, 48, 48, 48, 48, 48]

Candidates sorted by complexity= [0.07     0.43     0.57     0.18     0.82    
0.32     0.68     0.93
 0.02375  0.97625  0.47625  0.66375  0.77375  0.27375  0.72625  0.10125
 0.89875  0.39875  0.60125  0.14875  0.85125  0.22625  0.52375  0.695
 0.305    0.33625  0.83625  0.16375  0.58625  0.41375  0.91375  0.08625
 0.71125  0.28875  0.78875  0.21125  0.53875  0.46125  0.96125  0.03875
 0.055    0.945    0.445    0.555    0.195    0.805    0.35125  0.64875
 0.0075   0.3675   0.6325   0.9925   0.4925   0.5075   0.7575   0.2575
 0.2425   0.1175   0.8825   0.3825   0.6175   0.1325   0.8675   0.7425
 0.101875 0.148125 0.601875 0.398125 0.898125 0.726875 0.976875 0.773125
 0.226875 0.523125 0.476875 0.851875 0.273125 0.351875 0.789375 0.039375
 0.960625 0.460625 0.539375 0.210625 0.023125 0.289375 0.710625 0.085625
 0.914375 0.414375 0.585625 0.164375 0.835625 0.648125 0.679375 0.070625
 0.820625 0.335625 0.008125 0.991875 0.491875 0.508125 0.241875 0.758125
 0.258125 0.741875 0.116875 0.883125 0.383125 0.616875 0.133125 0.866875
 0.366875 0.633125 0.054375 0.945625 0.445625 0.554375 0.195625 0.804375
 0.304375 0.695625 0.929375 0.429375 0.570625 0.179375 0.320625 0.664375]

x_0 has rank r= 7  out of  128  candidates

\end{verbatim}

\section{Local density=stretch-based reconstruction }
\setcounter{figure}{0}

\begin{figure}
    \includegraphics[width=0.65\textwidth]{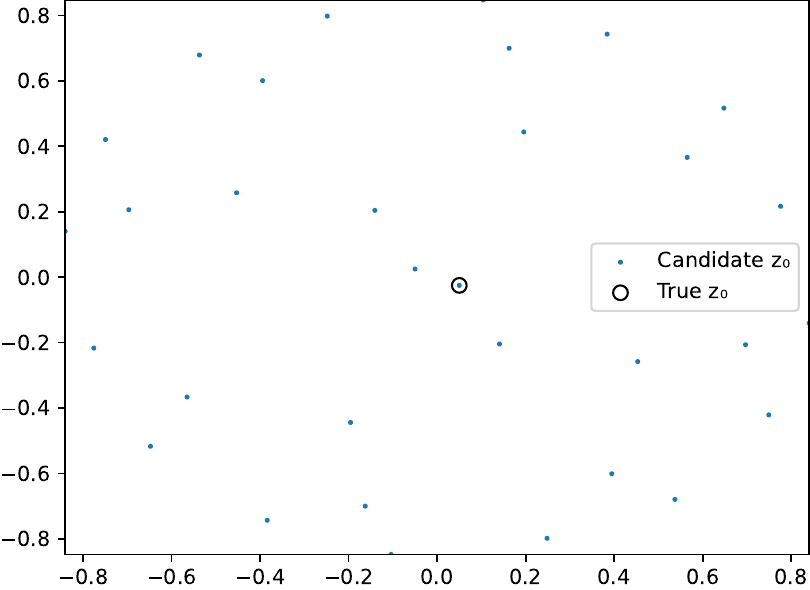}
     \includegraphics[width=0.65\textwidth]{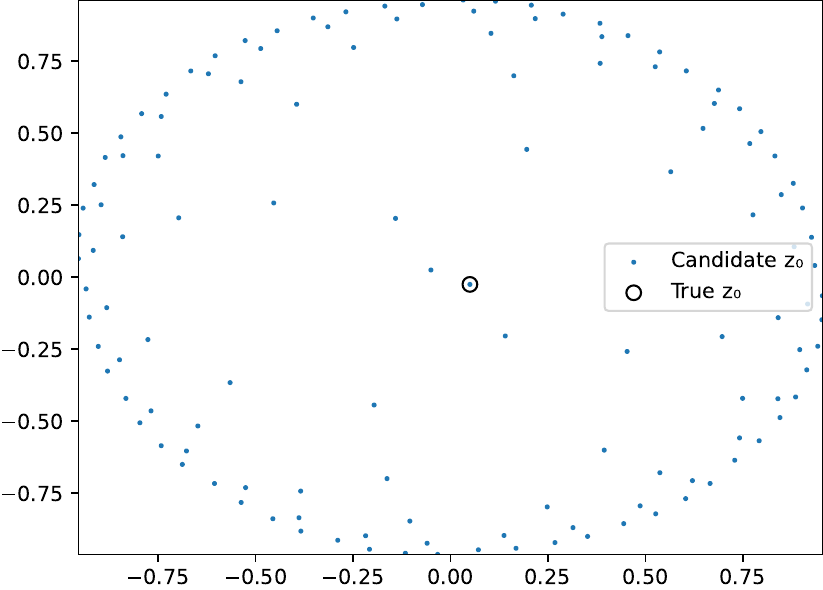}
  \caption{
  Density method with Julia map. For $n=5$ ($top$) the candidates are roughly equidistant from each other, so the density method does not work well. For $n=7$ ($bottom$), the true starting value is in a low-density region, and hence the density method works well. 
    }\label{fig:julia_n58}
\end{figure}

A density-based approach to reconstructing $x$ could be as follows:
If we ask what is the $\eps$-smeared probability $\rho(x)$ that $y$ maps to an $\eps^2$-neighborhood of $x$,
the answer for $\eps\to 0$ is $\rho(x)\propto 1/g'(x)$ for every $x\in\cX_0$,
where $g'$ can be computed from $f'$ via the chain rule.
Keeping $\eps$ small but finite gives a more smoothed and possibly more interpretable $\rho$.

Equivalently, in the forward direction, we could take a narrow Gaussian $x\sim\cN(x,\eps^2)$ 
around potential initial values $x$, which implies approximately $y\sim\cN(g(x),[g'(x)\eps]^2)$.
We call $g'(x)$ the stretch factor of $x$.
Choose $x\in\cX_0$ such that $y$ has minimal variance = smallest $g'(x)$ = largest ``density'' $\rho(x):=1/g'(x)$.
More precisely, sort $\cX_0$ w.r.t.\ increasing increasing $g'(x)$, then encode the rank $r$ of $x$ in $\log_2 r$ bits. Now let us consider how this method would apply to the maps we have studied in the main text.

\paradot{Logistic map}
For the logistic map with $n\to\infty$, 
one can show that $x\in\cX_0$ are neatly Beta(0.5,0.5) distributed,
such that $\rho(x)\propto\sqrt{x(1-x)}$,
so the density criterion is not useful. Hence this stretch method would not aid in  retrodiction for the logistic map with $\mu=4$. However, for the logistic map with $\mu<4$,  the point distribution $\cX_0$ and corresponding $g'(x)$ are more chaotic for most $\mu$, and the density criterion should yield shorter codes.

\paradot{Tent and Bernoulli maps}
The Tent and Bernoulli maps repeatedly and linearly stretch out an $\eps$-neighborhood of any $x$ by a factor of $2$,
so $|g'(x)|$ is constant independent of $n$, so this criterion is vacuous. Therefore, this stretch method will not apply to the Tent and Bernoulli maps and will not aid in  retrodiction. 

\paradot{Julia map}
We expect the complex Julia iteration $z\leftarrow z^2+c$ to get interesting results.
Since iteration is holomorphic, angles are preserved, hence $\rho(x)=1/|g'(z)|_2^2=|4z|^{-1}$.
In general, for $d>1$ we have $\rho(x)=1/|\text{det}(\nabla g(x))|$ by a similar line of arguments. Therefore, we expect this stretch method to aid in retrodiction of this map.

\end{document}